\documentclass{amsart}
\usepackage{fullpage}
\usepackage{epsfig, amsmath,xspace}
\usepackage{amssymb,amscd}
\usepackage[all,cmtip]{xy}
\usepackage{hyperref}

\newtheorem{thm}{Theorem}[section]
\newtheorem{lem}[thm]{Lemma}
\newtheorem{prop}[thm]{Proposition}
\newtheorem{cor}[thm]{Corollary}
\theoremstyle{definition}
\newtheorem{defn}[thm]{Definition}
\newtheorem{defnandnot}[thm]{Definition and Notation}

\newtheorem{ex}[thm]{Example}

\theoremstyle{remark}
\newtheorem{rem}[thm]{Remark}
\numberwithin{equation}{section}
\def\ie{\emph{i.e.}}

\newdir{ >}{{}*!/-7pt/@{>}}

\newcommand{\Sh}{{\mathrm{Sh}}}

\newcommand{\kmod}{k\text{-mod}}
\newcommand{\deltaepi}{{\Delta^{\mathrm{epi}}}}
\newcommand{\deltaepin}[1]{\Delta^{\mathrm{epi}}_{#1}}
\newcommand{\bepi}{b^\mathrm{epi}}
\newcommand{\epin}{\mathrm{Epi}_n}
\newcommand{\epi}{\mathrm{Epi}}

\def\L{\mathcal{L}}

\def\lra{\longrightarrow}
\def\id{\mathrm{id}}
\def\leq{\leqslant}
\def\geq{\geqslant}

\def\ra{\rightarrow}

\begin{document}
\title[$E_n$-homology as functor homology]{An interpretation of
  $E_n$-homology as functor homology}
\author{Muriel Livernet}
\address{Universit\'e Paris 13, CNRS, UMR 7539 LAGA, 99 avenue
  Jean-Baptiste Cl\'ement, 93430 
  Villetaneuse, France}
\email{livernet@math.univ-paris13.fr}
\author{Birgit Richter}
\address{Fachbereich Mathematik der Universit\"at Hamburg,
Bundesstra{\ss}e 55, 20146 Hamburg, Germany}
\email{richter@math.uni-hamburg.de}
\thanks{The first author thanks MIT and the Clay Institute for hosting
  her and  Haynes Miller for conversations on $E_n$-algebras. The
  second author thanks the  
Institut Galil\'ee of Universit\'e Paris 13 for an invitation as professeur 
invit\'e that led to this work. We are grateful to  Benoit Fresse for catching 
a serious sign error and to Aur\'elien Djament for suggesting a different 
setting for our proof of proposition 4.4. We thank the referee for
his/her careful  reading of the paper. }
\date{\today}
\begin{abstract}
We prove that $E_n$-homology of non-unital commutative algebras can be
described as functor homology
when one considers functors from a certain category of planar trees
with $n$ levels. For different $n$ these homology theories are
connected by natural maps, ranging from Hochschild homology  and its
higher order versions to Gamma homology.
\end{abstract}
\maketitle
\section{Introduction}
By neglect of structure, any commutative and associative algebra can
be considered as an associative algebra. More generally, we can
view such an algebra as an $E_n$-algebra, \ie, an algebra over an
operad in chain complexes that is weakly equivalent to the chain
complex of the little-$n$-cubes operad of \cite{BV} for $1 \leq n \leq
\infty$. Hochschild homology is a classical homology theory for
associative algebras and hence it can be applied to commutative
algebras as well. Less classically, Gamma homology  \cite{RoWh} is a homology
theory for $E_\infty$-algebras and Gamma homology of commutative
algebras plays an important role in the obstruction theory for
$E_\infty$ structures on ring spectra \cite{Ro,GH,BR} and its structural
properties are rather well understood \cite{RiRo}.

It is desirable to have a good understanding of the appropriate
homology theories in the intermediate range, \ie, for $1 < n < \infty$. A
definition of $E_n$-homology for augmented commutative algebras is due
to Benoit Fresse \cite{F} and the main topic of this paper is to prove
that these homology theories possess an interpretation in terms of
functor homology. We extend the range of $E_n$-homology to functors
from a suitable category $\epin$ to modules in such a way that it
coincides with Fresse's theory when we consider a functor that belongs
to an augmented commutative algebra and show in Theorem
\ref{thm:mainthm} that $E_n$-homology can
be described as functor homology, so that the homology groups are
certain Tor-groups.

As a warm-up we show in section \ref{sec:barhomology} that bar
homology of a non-unital algebra can be expressed in terms of functor
homology for functors from the category of order-preserving
surjections to $k$-modules. In section \ref{sec:epin} we introduce our
categories of epimorphisms, $\epin$, and their relationship to planar
trees with $n$-levels. We introduce a definition of $E_n$-homology for
functors from $\epin$ to $k$-modules that coincides with Benoit
Fresse's definition of $E_n$-homology of a  non-unital commutative
algebra, $\bar{A}$, when we apply our version of $E_n$-homology to
a suitable functor, $\L(\bar{A})$. We describe a spectral sequence
that has tensor products of bar homology groups as input and converges
to $E_2$-homology. Section \ref{sec:comp} is the technical heart of
the paper. Here we prove that $E_n$-homology has a Tor
interpretation. The proof of the acyclicity of a family of suitable
projective generators is an inductive argument that uses 
homology of small categories.

For varying $n$, the $E_n$-homology theories are related to each
other via a sequence of maps
$$H_*^{E_1} \rightarrow H_*^{E_2} \rightarrow H_*^{E_3} \rightarrow
\ldots$$ In a different context it is well known that the
stabilization map from Hochschild homology to Gamma homology can be
factored over so called higher order Hochschild homology \cite{P}:
for a commutative algebra $A$ there is a sequence of maps connecting
Hochschild homology of $A$, $HH_*(A)$, to Hochschild homology of
order $n$ of $A$ and finally to Gamma homology of $A$,
$H\Gamma_{*-1}(A)$. 
We explain how higher order Hochschild homology is related to
$E_n$-homology  for $n$ ranging from $1$ to $\infty$ in
\ref{subsec:comp}.

In the following we fix a commutative ring with unit, $k$. For a set
$S$ we denote by $k[S]$ the free $k$-module generated by $S$. If $S
= \{s\}$, then we write $k[s]$ instead of $k[\{s\}]$.
\section{Tor interpretation of bar homology}
\label{sec:barhomology}
We interpret the bar homology of a functor from the category of finite
sets and order-preserving surjections to the category of $k$-modules
as a $\mathrm{Tor}$-functor.

For unital $k$-algebras,
the complex for the Hochschild homology of the algebra can be viewed
as the chain complex associated to a simplicial object. In the absense
of units, this is no longer possible.

Let $\bar{A}$ be a non-unital $k$-algebra. The bar-homology of
$\bar{A}$, $H^{\mathrm{bar}}_*(\bar{A})$,  is defined as the homology
of the  complex
$$ C^{\mathrm{bar}}_*(\bar{A}):  \ldots \rightarrow \bar{A}^{\otimes n+1}
\stackrel{b'}{\longrightarrow}  \bar{A}^{\otimes n}
\stackrel{b'}{\longrightarrow} \ldots \stackrel{b'}{\longrightarrow}
\bar{A} \otimes \bar{A} \stackrel{b'}{\longrightarrow} \bar{A}$$
with $C^{\mathrm{bar}}_n(\bar{A}) = \bar{A}^{\otimes n+1}$ and $b' =
\sum_{i=0}^{n-1} (-1)^i d_i$ where $d_i$ applied to $a_0 \otimes
\ldots \otimes a_n \in \bar{A}^{\otimes n+1}$ is $a_0 \otimes \ldots
\otimes a_ia_{i+1} \otimes  \ldots \otimes a_n$.
\medskip

The category of non-unital associative $k$-algebras is equivalent to
the  category of augmented $k$-algebras.
If one replaces $\bar{A}$ by $A=\bar{A}\oplus k$, then
$C^{\mathrm{bar}}_n(\bar{A})$ corresponds to the reduced Hochschild
complex of  $A$ with coefficients in the trivial module $k$, shifted
by one:  $H^{\mathrm{bar}}_*(\bar A)=HH_{*+1}(A,k)$, for $*\geq 0$.

\begin{defn}
Let $\deltaepi$ be the category whose objects are the sets $[n] =
\{0,\ldots,n\}$ for $n \geq 0$ with the ordering $0 < 1 < \ldots <
n$ and whose morphisms are order-preserving surjective functions. We
will call covariant functors $F\colon \deltaepi \rightarrow  \kmod$
\emph{$\deltaepi$-modules}.
\end{defn}
We have the basic order-preserving surjections $d_i\colon [n]
\rightarrow [n-1], 0 \leq i \leq n-1$ that are given by
$$ d_i(j) = \left\{\begin{array}{rl}
j, & j \leq i, \\
j-1, & j > i.
\end{array} \right.$$
Any order-preserving surjection is a composition of these basic ones.
The generating morphisms $d_i$ in $\deltaepi$ correspond to the face
maps in the standard simplicial model of the $1$-sphere with the
exception of the last face map. The standard simplicial identities imply
that  $b' = \sum_{i=0}^{n-1} (-1)^i F(d_i)$ satisfies $b'^2=0$. It
justifies the following definition. 

\begin{defn}\label{def:barhom}
We define the \emph{bar-homology of a $\deltaepi$-module $F$}  as the
homology of the complex
$C^{\mathrm{bar}}_*(F)$ with $C^{\mathrm{bar}}_n(F) = F[n]$ and
    differential $b' = \sum_{i=0}^{n-1} (-1)^i F(d_i)$.
\end{defn}

For a non-unital algebra $\bar{A}$ the functor $\L(\bar{A})$ that
assigns $\bar{A}^{\otimes (n+1)}$ to $[n]$ and $\L(d_i)(a_0 \otimes
\ldots \otimes a_n) = a_0 \otimes \ldots \otimes a_ia_{i+1} \otimes
\ldots \otimes a_n$ ($0 \leq i \leq n-1$) is a $\deltaepi$-module.
In that case, $C^{\mathrm{bar}}_*(\L(\bar{A})) =
C^{\mathrm{bar}}_*(\bar{A})$.

In the following we use the machinery of functor homology as in
\cite{PR}. Note that the category of $\deltaepi$-modules has enough
projectives: the representable functors $(\deltaepi)^n \colon \deltaepi
\rightarrow \kmod$ with $(\deltaepi)^n[m] = k[\deltaepi([n],[m])]$ are
easily seen to be projective objects and each $\deltaepi$-module
receives a surjection from a sum of representables. The analogous
statement is true for contravariant functors from $\deltaepi$ to the
category of $k$-modules where we can use the functors $\deltaepin{n}$
with  $\deltaepin{n}[m] = k[\deltaepi([m],[n])]$ as
projective objects.

We call the cokernel of the map between contravariant representables
$$ (d_0)_*\colon \deltaepin{1} \rightarrow \deltaepin{0}$$
$\bepi$.
Note that $\deltaepin{0}[n]$ is free of rank one for all $n \geq 0$
because there is just one map in $\deltaepi$ from $[n]$ to $[0]$ for
all $n$. Furthermore, $\deltaepin{1}[0]$ is the zero module, because
$[0]$ cannot surject onto $[1]$. Therefore
$$ \bepi[n] \cong \left\{
  \begin{array}{cr}
    0 & \text{ for } n > 0, \\
 k & \text{ for } n = 0.
  \end{array} \right. $$
\begin{prop}\label{prop:n=1}
For any $\deltaepi$-module $F$
\begin{equation}
H_p^{\mathrm{bar}}(F) \cong \mathrm{Tor}^{\deltaepi}_p(\bepi, F)
\text{ for all } p \geq 0.
\end{equation}
\end{prop}
For the proof recall that a sequence of $\deltaepi$-modules and
natural transformations
\begin{equation} \label{eq:ses}
0 \ra F' \stackrel{\phi}{\longrightarrow} F
\stackrel{\psi}{\longrightarrow} F'' \rightarrow 0
\end{equation}
is \emph{short exact} if it gives rise to a short exact sequence of
$k$-modules
$$ 0 \ra F'[n] \stackrel{\phi[n]}{\longrightarrow} F[n]
\stackrel{\psi[n]}{\longrightarrow} F''[n] \rightarrow 0$$
for every $n \geq 0$.
\begin{proof}
We have to show that $H_*^{\mathrm{bar}}(-)$ maps short exact sequences of
$\deltaepi$-modules to long exact sequences, that $H_*^{\mathrm{bar}}(-)$
vanishes on projectives in  positive degrees and that $H_0^{\mathrm{bar}}(F)$ and
$\bepi \otimes_{\deltaepi} F$ agree for all $\deltaepi$-modules $F$.

A short exact sequence as in \eqref{eq:ses} is sent to a short exact
sequence of chain complexes
$$ \xymatrix@1{
{0} \ar[r] & {C_*^{\mathrm{bar}}(F')}
\ar[rr]^{C_*^{\mathrm{bar}}(\phi)} & & {C_*^{\mathrm{bar}}(F)}
\ar[rr]^{C_*^{\mathrm{bar}}(\psi)} & & {C_*^{\mathrm{bar}}(F'')} \ar[r] &
{0} }$$
and therefore the first claim is true.

In order to show that $H_*^{\mathrm{bar}}(P)$ is trivial in positive degrees
for any projective $\deltaepi$-module $P$ it suffices to show that the
representables $(\deltaepi)^n$ are acyclic. In order to prove this claim
we construct an explicit chain homotopy.

Let $f \in (\deltaepi)^n[m]$ be a generator, \ie, a surjective
order-preserving map from $[n]$ to $[m]$. Note that $f(0) = 0$. We
can codify such a map by its fibres, \ie, by an $(m+1)$-tuple of
pairwise disjoint  subsets  $(A_0, \ldots, A_m)$ with $A_i \subset
[n]$, $0 \in A_0$ and $\bigcup_{i=0}^{m} A_i = [n]$ such that $x <y$
for $x \in A_i$ and $y \in A_j$ with $i<j$. With this notation
$d_i(A_0,\ldots,A_m)=(A_0,\ldots,A_{i-1},A_i\cup
A_{i+1},\ldots,A_m)$.

We define the chain homotopy $h\colon \deltaepi([n],[m]) \rightarrow
\deltaepi([n],[m+1])$ as
\begin{equation} \label{eq:htp}
h(A_0,\ldots,A_m) := \left\{ \begin{array}{cl}
0 & \text{ if } A_0 = \{0\}, \\
(0,A_0',A_1,\ldots, A_m) & \text{ if } A_0 = \{0\} \cup A_0', A_0'
\neq \varnothing.
\end{array} \right.
\end{equation}
If $A_0 = \{0\}$, then
$$ (b'\circ h + h \circ b')(\{0\},\ldots,A_m) = 0 + h \circ
b'(\{0\},\ldots,A_m) = h(\{0\} \cup A_1, \ldots, A_m) =
(\{0\},\ldots,A_m).$$ In the other case a direct calculation shows
that $ (b'\circ h + h \circ b')(A_0,\ldots,A_m) =
\id(A_0,\ldots,A_m)$.

It remains to show that both homology theories coincide in degree
zero. By definition $H_0^{\mathrm{bar}}(F)$ is the cokernel of the map
$$ F(d_0)\colon F[1] \longrightarrow F[0].$$
A Yoneda-argument \cite[17.7.2(a)]{Sch} shows that the tensor product
$\deltaepin{n} \otimes_{\deltaepi} F$ is naturally
isomorphic to $F[n]$ and hence the above cokernel is the cokernel
of the map
$$ ((d_0)_* \otimes_{\deltaepi} \id) \colon \deltaepin{1}
\otimes_{\deltaepi} F \longrightarrow \deltaepin{0} \otimes_{\deltaepi}
F.$$
As tensor products are right-exact \cite[17.7.2 (d)]{Sch}, the cokernel of the
above map is isomorphic to
$$ \mathrm{coker}((d_0)_* \colon \deltaepin{1} \rightarrow
\deltaepin{0}) \otimes_{\deltaepi}
F = \bepi \otimes_{\deltaepi}
F = \mathrm{Tor}_0^\deltaepi(\bepi, F).$$
\end{proof}

\begin{rem}\label{R:gradedn=1} Note that the previous proof applies in
  a graded  context. More precisely,
assume every element $i \in  [n]$ comes with a grading $d(i) \in
\mathbb{N}_0$. In the sequel we will consider surjective maps
$\phi:X\rightarrow [n]$ where $X$ is a graded set and
$d(i)=\sum_{x|\phi(x)=i} d(x)$.
The degree of a subset $A$ of $[n]$  is $d(A)=\sum_{i\in A} d(i)$. To
this grading one can associate the complex
$(C_m^{X,\phi}((\deltaepi)^n)=\bigoplus_{f\in\deltaepi([n],[m])}
k[f],b')$, with $d_i(A_0,\ldots,A_m)=(-1)^{\sum_{j=0}^i
d(A_j)}(A_0,\ldots,A_i\cup A_{i+1},\ldots,A_m)$ and
$b'=\sum_{i=0}^{m-1} (-1)^i d_i$. As the $d_i$'s still satisfy the
simplicial identities, $b'$ is of square $0$. The complexes
$C_*^{X,\phi}((\deltaepi)^n)$ are acyclic, as one can see using  the homotopy
\begin{equation*}
h(A_0,\ldots,A_m) := \left\{ \begin{array}{cl}
0 & \text{ if } A_0 = \{0\}, \\
(-1)^{d(0)}(0,A_0',A_1,\ldots, A_m) & \text{ if } A_0 = \{0\} \cup
A_0', A_0' \neq \varnothing.
\end{array} \right.
\end{equation*}
\end{rem}

\section{Epimorphisms and trees} \label{sec:epin}
Planar level trees are used in \cite{B}, \cite{F} and
\cite[3.15]{Be} as a means to codify $E_n$-structures. An $n$-level
tree is a planar level tree with $n$ levels. We will use categories
of planar level trees in order to gain a description of
$E_n$-homology as functor homology. If $\mathcal{C}$ is a small
category, then we denote the nerve of $\mathcal{C}$ by
$N\mathcal{C}$.

\begin{defn}
Let $n \geq 1$ be a natural number. The category $\epin$ has as
objects the elements of $N_{n-1}(\deltaepi)$, \ie, sequences
\begin{equation} \label{eq:objects}
\xymatrix@1{{[r_n]} \ar[r]^{f_n} & {[r_{n-1}]} \ar[r]^{f_{n-1}} &
  {\ldots} \ar[r]^{f_2} & {[r_1]}}
\end{equation}
with $[r_i] \in \deltaepi$ and surjective order-preserving  maps $f_i$. A
morphism in $\epin$ from the above object to an object
$\xymatrix@1{{[r'_n]} \ar[r]^{f'_n} & {[r'_{n-1}]} \ar[r]^{f'_{n-1}} &
  {\ldots} \ar[r]^{f'_2} & {[r'_1]}} $
consists of surjective maps $\sigma_i\colon [r_i] \rightarrow [r'_i]$
for $1 \leq i \leq n$ such that  $\sigma_1\in\deltaepi$ and for all $2
\leq i \leq n$  the map  $\sigma_i$
is  order-preserving on the fibres $f_i^{-1}(j)$ for all
$j \in [r_{i-1}]$  and such that the
diagram
$$ \xymatrix{{[r_n]} \ar[d]^{\sigma_n} \ar[r]^{f_n} & {[r_{n-1}]}
  \ar[d]^{\sigma_{n-1}} \ar[r]^{f_{n-1}} &
  {\ldots}  \ar[r]^{f_2} & {[r_1]} \ar[d]^{\sigma_{1}} \\
{[r'_n]} \ar[r]^{f'_n} & {[r'_{n-1}]} \ar[r]^{f'_{n-1}} &
  {\ldots} \ar[r]^{f'_2} & {[r'_1]}}
$$
commutes.
\end{defn}
As an example, consider the object $\xymatrix@1{{[2]} \ar[r]^{\id} &
  {[2]}}$ in $\mathrm{Epi}_2$ which can be viewed as the 2-level tree

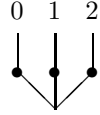
\vspace{1.3cm}
\begin{center}
\begin{picture}(2,3)
\setlength{\unitlength}{1cm}
\put(1,0){\line(-1,1){0.5}}
\put(1,0){\line(1,1){0.5}}
\put(1,0){\line(0,1){1}}
\put(0.5,0.5){\line(0,1){0.5}}
\put(1.5,0.5){\line(0,1){0.5}}
\put(1.4,1.2){2}
\put(0.9,1.2){1}
\put(0.4,1.2){0}
\put(0.4,0.4){$\bullet$}
\put(0.9,0.4){$\bullet$}
\put(1.4,0.4){$\bullet$}
\end{picture}
\end{center}
Possible maps from this object to $\xymatrix@1{{[2]} \ar[r]^{d_0} &
  {[1]}}$ are \,
$\vcenter{\xymatrix{
{[2]} \ar[r]^{\id} \ar[d]_{\id}&   {[2]} \ar[d]^{d_0} \\
{[2]} \ar[r]^{d_0} &  {[1]}}}$ \,  and \, $\vcenter{\xymatrix{
{[2]} \ar[r]^{\id} \ar[d]_{(0,1)}&   {[2]} \ar[d]^{d_0} \\
{[2]} \ar[r]^{d_0} &  {[1]}}}$ \, where $(0,1)$ denotes the transposition
that permutes $0$ and $1$. For $\sigma_1 = d_1$ there is no possible
$\sigma_2$ to fill in the diagram.

If $n=1$, then $\mathrm{Epi}_1$ coincides with the category
$\deltaepi$. Note that there is a functor $\iota_n\colon \deltaepi =
\mathrm{Epi}_1 \rightarrow \epin$ for all $n \geq 1$ with
$$\iota_n([m]) := \xymatrix@1{{[m]} \ar[r] & {[0]} \ar[r] & {\ldots}
  \ar[r] & {[0].}}$$
We call trees of the form $\iota_n([m])$ \emph{palm trees with $m+1$
  leaves.} More generally we have functors connecting the various
categories of planar level trees.

\begin{lem} \label{lem:epin-1toepin}
For all $n>k\geq 1$  there are functors $\iota_n^k \colon
\mathrm{Epi}_k \rightarrow \epin$, with
$$\iota_n^k (\xymatrix@1{[r_k] \ar[r]^{f_k} & \ldots \ar[r]^{f_2} &
  [r_1]}) = \xymatrix@1{[r_k] \ar[r]^{f_k} & \ldots \ar[r]^{f_2} &
  [r_1] \ar[r] & [0] \ar[r] & \ldots \ar[r] & [0]}$$
on objects, with the canonical extension to morphisms. \qed
\end{lem}
\begin{rem} \label{rem:fork} The maps $\iota_n^k$ correspond to
  iterated suspension morphisms in
\cite[4.1]{B}. There is a different way of mapping a planar tree
with $n$ levels to one with $n+1$ levels, by sending
$\xymatrix@1{[r_n] \ar[r]^{f_n} & \ldots \ar[r]^{f_2} &
  [r_1]}$ to $\xymatrix@1{[r_n] \ar[r]^{\id_{[r_n]}} & [r_n]
  \ar[r]^{f_n} &  \ldots \ar[r]^{f_2} &
  [r_1]}$.
\end{rem}
\begin{defn}
We call trees of the form $\xymatrix@1{[r_n] \ar[r]^{\id_{[r_n]}} &
[r_n]   \ar[r]^{f_n} &  \ldots \ar[r]^{f_2} &   [r_1]}$ \emph{fork
trees}.
\end{defn}
Fork trees will  need special attention later when we prove that
representable functors are acyclic.

For any $\Sigma_*$-cofibrant operad $\mathcal P$ there exists a
homology  theory for
$\mathcal P$-algebras which is denoted by $H_*^{\mathcal P}$ and is
called $\mathcal P$-homology. Fresse studies the particular case of
$\mathcal P=E_n$ a differential graded operad quasi-isomorphic to the
chain operad of the little $n$-disks operad. He proves that for any
commutative algebra the $E_n$-homology coincides with the homology of
its  $n$-fold bar construction. In fact, his result is more general
since he defines an analogous $n$-fold  bar construction for
$E_n$-algebras and proves the result for any $E_n$-algebra
in \cite[theorem 7.26]{F}.

We consider the $n$-fold  bar construction of a non-unital
commutative $k$-algebra $\bar{A}$, $B^n(\bar{A})$, as an $n$-complex, such that
$$ B^n(\bar{A})_{(r_n,\ldots,r_1)} =
\bigoplus_{[r_n] \stackrel{f_n}{\rightarrow} \ldots
\stackrel{f_2}{\rightarrow} [r_1] \in \epin} \bar{A}^{\otimes
(r_n+1)}.$$ The differential in $B^n(\bar{A})$ is the total
differential associated to $n$-differentials $\partial_1, \ldots,
\partial_n$ such that $\partial_n$ is built out of the
multiplication in $\bar{A}$, $\partial_{n-1}$ corresponds to the
shuffle multiplication on $B(\bar{A})$ and so on. We describe the
precise setting in a slightly more general context.

In order to extend the $\mathrm{Tor}$-interpretation of bar homology
of $\deltaepi$-modules to functors from $\epin$ to modules (alias
$\epin$-modules) we describe the $n$ kinds of face maps  for $\epin$
in detail by considering diagrams of the form
\begin{equation}  \label{eq:di}
\xymatrix{
{[r_n]} \ar[r]^{f_n} \ar[d]_{\tau_{n}^{i,j}}& {[r_{n-1}]} \ar[r]^{f_{n-1}}
  \ar[d]_{\tau_{n-1}^{i,j}}&
  {\ldots}\ar[r]^{f_{j+2}}&[r_{j+1}] \ar[d]_{\tau_{j+1}^{i,j}}
  \ar[r]^{f_{j+1}} & [r_j]  \ar[r]^{f_{j}} \ar[d]_{d_i} & 
{[r_{j-1}]} \ar[d]_{\id}
  \ar[r]^{f_{j-1}} & {\ldots} \ar[r]^{f_2}  &  {[r_1]} \ar[d]_{\id}\\
{[r_n]} \ar[r]^{g_n} & {[r_{n-1}]} \ar[r]^{g_{n-1}} &
  {\ldots}\ar[r]^{g_{j+2}}& [r_{j+1}]\ar[r]^{g_{j+1}} & [r_j-1]
  \ar[r]^{g_{j}} & {[r_{j-1}]} \ar[r]^{f_{j-1}} &{\ldots}
  \ar[r]^{f_2} &  {[r_1].}
}
\end{equation}
Given the object in the first row, it is not always possible to
extend $(d_i\colon [r_j]\rightarrow [r_j-1], \id_{[r_{j-1}]},\ldots,
\id_{[r_1]})$ to a  morphism in $\epin$: we have to find
order-preserving surjective maps $g_k$  for $j \leq k \leq n$ and
bijections $\tau_{k}^{i,j}\colon [r_k] \rightarrow [r_k]$ that are
order-preserving on the  fibres of $f_k$ for $j+1 \leq k \leq n$
such that the diagram commutes.

By convention we denote the constant map $[r_1]\rightarrow [0]$ by $f_1$.

\begin{lem}\label{lemma:di}
\begin{enumerate}
\item[]
\item
There is a unique order-preserving surjection
$g_j\colon [r_j-1] \rightarrow [r_{j-1}]$ with $g_j \circ
d_i = f_j$ if and only if $f_j(i) = f_j(i+1)$. When it exists, $g_j$
is denoted by $f_j|_{i=i+1}$. 
\item
If $f_j(i)=f_j(i+1)$ then we can extend the diagram to one of the form
\eqref{eq:di} so that $\tau_{j+1}^{i,j}$ is a shuffle of the fibres
$f_{j+1}^{-1}(i)$ and $f_{j+1}^{-1}(i+1)$. Each choice of a
$\tau_{j+1}^{i,j}$ uniquely determines the maps $\tau_{k}^{i,j}$ for
all $j+1 < k \leq n$.
\item If $f_j(i)=f_j(i+1)$ then  each choice of a
$\tau_{j+1}^{i,j}$ uniquely determines the maps $g_k$ for $k\geq j$.
The diagram  \eqref{eq:di} takes the following form
\begin{equation}  \label{eq:di2}
\xymatrix{
{[r_n]} \ar[r]^{f_n} \ar[d]_{\tau_{n}^{i,j}}& {[r_{n-1}]} \ar[r]^{f_{n-1}}
  \ar[d]_{\tau_{n-1}^{i,j}}&
  {\ldots}\ar[r]^{f_{j+2}}&[r_{j+1}] \ar[d]_{\tau_{j+1}^{i,j}}
  \ar[r]^{f_{j+1}} & [r_j]  \ar[r]^{f_{j}} \ar[d]_{d_i} & 
{[r_{j-1}]} \ar[d]_{\id}
  \ar[r]^{f_{j-1}} & {\ldots} \ar[r]^{f_2}  &  {[r_1]} \ar[d]_{\id}\\
{[r_n]} \ar[r]^{g_n^\tau} & {[r_{n-1}]} \ar[r]^{g_{n-1}^\tau} &
  {\ldots}\ar[r]^{g_{j+2}^\tau}& [r_{j+1}]\ar[r]^{d_if_{j+1}} &
  [r_j-1]  \ar[r]^{f_{j}|_{i=i+1}} & {[r_{j-1}]} 
\ar[r]^{f_{j-1}} &{\ldots}
  \ar[r]^{f_2} &  {[r_1].}
}
\end{equation}

  \end{enumerate}
\end{lem}

\begin{proof}
If there is such a map $g_j$, then
$f_j(i+1) = g_j \circ d_i(i+1) = g_j \circ d_i(i) = f_j(i)$. As $f_j$
is order-preserving, it is determined by the cardinalities of its
fibres. The decomposition of morphisms in the simplicial category then
ensures that we can factor $f_j$ in the desired way.

For the third claim, assume that $g_j$ exists with the properties
mentioned  in (a).
As $g_{j+1}$ and $d_i \circ f_{j+1}$ are both order-preserving maps
from $[r_{j+1}]$ to $[r_j-1]$, they are determined by the
cardinalities of the fibres and thus they have to agree.
Then $\tau_{j+1}^{i,j} = \id_{[r_{j+1}]}$
extends the diagram up to layer $j+1$. For the higher layers we then
have to choose $g_k = f_k$ and $\tau_{k}^{i,j} = \id_{[r_k]}$.

In general, $\tau_{j+1}^{i,j}$ has to satisfy the conditions that it
is order-preserving on the fibres of $f_{j+1}$. If $A_i =
f_{j+1}^{-1}(i)$ then this implies that $\tau_{j+1}^{i,j}$ is an
$(A_0,\ldots, A_{r_j})$-shuffle. Furthermore we have that
$$ (d_i \circ f_{j+1})^{-1}(k) =
\left\{ \begin{array}{cl}
  A_k & \text{ if } k < i, \\
A_i \cup A_{i+1} & \text{ if } k = i, \\
A_{k+1} & \text{ if } k > i.
\end{array} \right. $$
Therefore $\tau_{j+1}^{i,j}$ has to map $A_0, \ldots, A_{i-1},
A_{i+2}, \ldots, A_{r_j}$ identically and is hence an
$(A_i,A_{i+1})$-shuffle.

If we fix a shuffle $\tau_{j+1}^{i,j}$, then the next permutation
$\tau_{j+2}^{i,j}$ has to be order-preserving on the fibres of
$f_{j+2}$, thus it is at most a shuffle of the fibres. In addition, it
has to  satisfy
\begin{equation} \label{eq:cond}
g_{j+2} \circ \tau_{j+2}^{i,j} = \tau_{j+1}^{i,j} \circ f_{j+2}.
\end{equation}
Again, as $g_{j+2}$ is order-preserving we have no choice but to take
the order-preserving map satisfying $|g_{j+2}^{-1}(k)| =
|(\tau_{j+1}^{i,j} \circ f_{j+2})^{-1}(k)|$, for all $k\in [r_{j+1}]$.
By \eqref{eq:cond} we know that
$\tau_{j+2}^{i,j}$ has to send $f_{j+2}^{-1}(k)$ to
$g_{j+2}^{-1}(\tau_{j+1}^{i,j}(k))$ and this determines
$\tau_{j+2}^{i,j}$. A proof by induction shows the general claim  in
(b).
\end{proof}

In the following we will extend the notion of $E_n$-homology for
commutative non-unital $k$-algebras to $\epin$-modules. Again thanks to
Fresse's theorem \cite[theorem 7.26]{F},
the $E_n$-homology and the homology of the $n$-fold bar construction
of a commutative algebra coincide.

Definitions \ref{def:endiff} and \ref{def:enchains} describe a
complex while lemma \ref{lem:chain1} asserts the complex property.

\begin{defnandnot} \label{def:endiff}
Let $t:[r_n] \stackrel{f_n}{\rightarrow} \ldots
\stackrel{f_{j+1}}{\rightarrow} [r_j] \stackrel{f_j}{\rightarrow}
\ldots \stackrel{f_2}{\rightarrow} [r_1]$ be an $n$-level tree.

The \emph{degree of $t$}, denoted by $d(t)$, is the number of its edges, that
is $\sum_{i=1}^n (r_i+1)$.

For fixed $1\leq j\leq n$ and $i\in[r_j]$ let $t_{j,i}$ be the
$(n-j)$-level tree defined by the $j$-fibre of $t$ over $i\in [r_j]$:
$$(f_{j+1}f_{j+2}\ldots f_{n})^{-1}(i) \stackrel{f_n}{\rightarrow}
\ldots  \stackrel{f_{k+1}}{\rightarrow} (f_{j+1}\ldots
f_{k})^{-1}(i) \stackrel{f_k}{\rightarrow} \ldots
\stackrel{f_{j+2}}{\rightarrow} f_{j+1}^{-1}(i).$$ Conversely a tree
$t$ can be recovered from its $1$-fibres, that is
$t=[t_{1,0},\ldots,t_{1,r_1}]$ and $d(t)=\sum_{i=0}^{r_1}
(d(t_{1,i})+1) = r_1 + 1 + \sum_{i=0}^{r_1} d(t_{1,i})$.

Let $F$ be an $\epin$-module. For fixed $1\leq j\leq n$ and $i\in[r_j]$
such that $f_j(i) = f_j(i+1)$ or
$j=1$  we define
$$d_i^j \colon F([r_n] \stackrel{f_n}{\rightarrow} \ldots
\stackrel{f_{j+1}}{\rightarrow} [r_j] \stackrel{f_j}{\rightarrow}
\ldots \stackrel{f_2}{\rightarrow} [r_1]) \longrightarrow
\bigoplus_{t'=[r_n] \stackrel{g_n}{\rightarrow}\ldots
\stackrel{g_{j+2}}{\rightarrow} [r_{j+1}]
\stackrel{d_if_{j+1}}{\rightarrow} [r_j-1]
\stackrel{f_j|_{i=i+1}}{\rightarrow} \ldots
\stackrel{f_2}{\rightarrow} [r_1] \in \epin} F(t')$$ as
\begin{equation} \label{eq:dij}
\begin{cases}
d_i^j = \sum_{\tau_{j+1}^{i,j} \in
\mathrm{Sh}(f_{j+1}^{-1}(i),f_{j+1}^{-1}(i+1))}
\varepsilon(\tau_{j+1}^{i,j};t_{j,i},t_{j,i+1}) 
F(\tau_n^{i,j},\ldots, \tau_{j+1}^{i,j}, d_i,\id, \ldots, \id)& \text{
  if }  j<n, \\
d_i^n=F(d_i,\id, \ldots, \id) & \text{if } j=n,
\end{cases}
\end{equation}
where the sign $\epsilon(\tau_{j+1}^{i,j};t_{j,i};t_{j,i+1})$ is
defined as follows: one expresses the $(n-j)$-level trees $t_{j,i}$
and $t_{j,i+1}$ as a sequence of $(n-j-1)$-level trees
$t_{j,i}=[t_1,\ldots,t_p]$ and $t_{j,i+1}=[t_{p+1},\ldots t_{p+q}]$;
the shuffle $\sigma=\tau_{j+1}^{i,j}$ is indeed a $(p,q)$-shuffle
and acts on $t$ by replacing the fibres $t_{j,i}$ and $t_{j,i+1}$ by
the fibre
$u_{j,i}=[t_{\sigma^{-1}(1)},\ldots,t_{\sigma^{-1}(p+q)}]$. The sign
$\varepsilon(\sigma;[t_1,\ldots,t_p], [t_{p+1},\ldots, t_{p+q}])$
picks up a factor of $(-1)^{(d(t_a)+1)(d(t_b)+1)}$ whenever
$\sigma(a) >\sigma(b)$ but $a<b$.

\end{defnandnot}

In \ref{ex:tree14} we treat an example with $n=3$.

\begin{defn} \label{def:enchains}
  \begin{itemize}
\item[]
\item
If $F$ is an $\epin$-module, then the \emph{$E_n$-chain complex of
$F$} is the $n$-fold chain complex whose $(r_n,\ldots,r_1)$ spot is
\begin{equation}
  \label{eq:cen}
  C^{E_n}_{(r_n,\ldots,r_1)}(F) =
\bigoplus_{[r_n] \stackrel{f_n}{\rightarrow} \ldots
\stackrel{f_2}{\rightarrow} [r_1] \in \epin} F([r_n]
\stackrel{f_n}{\rightarrow}  \ldots \stackrel{f_2}{\rightarrow}
[r_1]).
\end{equation}
The differential in the $j$-th coordinate is
$$ \partial_j \colon C^{E_n}_{(r_n,\ldots,r_j, \ldots, r_1)}(F)
\rightarrow C^{E_n}_{(r_n,\ldots,r_j-1, \ldots, r_1)}(F)$$ with
$$ \partial_j := \sum_{i|f_j(i) = f_j(i+1)} (-1)^{s_{j,i}} d_i^j,$$
 where $s_{j,i}$ is obtained as follows: drawing the tree $t$ on a
 plane with  its  root at the bottom, one can label its edges
-- from $1$ to $d(t)$ -- from bottom to top and left to right; the
integer $s_{j,i}$ is the label of  the right most top edge of the
tree $t_{j,i}$. For $j=n$ we use the convention that $s_{n,i}$ is
the label of the $i$-th leaf of $t$ for $0 \leq i \leq r_n$.
\item
The \emph{$E_n$-homology of $F$}, $H_*^{E_n}(F)$ is defined to be the
homology of the total complex associated to \eqref{eq:cen}.
\end{itemize}
\end{defn}
Note that for $n=1$ we recover definition \ref{def:barhom}, but with
$\partial_1=-b'$.
\begin{lem} \label{lem:chain1}
The $k$-modules $C^{E_n}_{(r_n,\ldots,r_1)}(F)$ constitute an
$n$-fold chain complex.
\end{lem}
We postpone the proof of the lemma until after example \ref{ex:tree14}.

\medskip

In order to prove the main theorem \ref{thm:mainthm}, we need
categories of $n$-trees depending on a fixed finite ordered set $X$
of graded elements denoted $\epi_n^X$. For any $x\in X$, $d(x) \in
\mathbb{N}_0$ will denote its degree. For any subset $A$ of $X$ the
degree $d(A)$ is the sum of the degrees of the elements of $A$, thus
for instance $d(X)=\sum_{x\in X} d(x)$. If $A,B$ is a pair of
disjoint subsets of $X$, then  one defines
$\epsilon(A;B)=\prod_{a\in A;b\in B; a>b} (-1)^{d(a)d(b)}$. One has
\begin{equation} \label{eq:epsilonAB}
\epsilon(A;B)\epsilon(B;A)=(-1)^{d(A)d(B)}.
\end{equation}

An object in the category $\epi_n^X$ is an $n$-level tree $t$
together with a surjection $\phi: X\rightarrow [r_n]$. Any such
element is denoted by $(t,\phi)$ and is called an $(X,n)$-level
tree. A morphism from $(t,\phi)$ to $(t',\phi')$ is a morphism
$\sigma\colon t\rightarrow t'$ in the category $\epi_n$ satisfying
$\phi'=\sigma_n\phi$. The following should be considered as a graded
version of \ref{def:endiff} and \ref{def:enchains}, the consistency
of the definition is stated in lemma \ref{lem:chain2}.

\begin{defn}  \label{def:enchainsX}
Let $(t,\phi): X \stackrel{\phi}{\rightarrow} [r_n]
\stackrel{f_n}{\rightarrow} \ldots \stackrel{f_{j+1}}{\rightarrow}
[r_j] \stackrel{f_j}{\rightarrow} \ldots \stackrel{f_2}{\rightarrow}
[r_1]$ be an $(X,n)$-level tree in $\epi^X_n$.

The \emph{degree of $(t,\phi)$}, denoted by $d(t)$, is the sum of the number
of its edges  and the degrees of elements of $X$,
$$ d(t) = \sum_{i=1}^n (r_i+1)+d(X).$$

For a fixed $1\leq j\leq n-1$ and $i\in[r_j]$ let $t_{j,i}$ be the
$(X_{j,i},n-j)$-level tree defined by the $j$-fibre of $t$ over $i\in [r_j]$:
$$X_{j,i}=(f_{j+1}f_{j+2}\ldots
f_{n}\phi)^{-1}(i)\stackrel{\phi}{\rightarrow}(f_{j+1}f_{j+2}\ldots
f_{n})^{-1}(i) \stackrel{f_n}{\rightarrow} \ldots
\stackrel{f_{k+1}}{\rightarrow} (f_{j+1}\ldots f_{k})^{-1}(i)
\stackrel{f_k}{\rightarrow} \ldots \stackrel{f_{j+2}}{\rightarrow}
f_{j+1}^{-1}(i).$$ For $j=n$, $t_{n,i}$ coincides with
$X_{n,i}=\phi^{-1}(i)$.

Conversely an $(X,n)$-level tree $(t,\phi)$ can be recovered from
its $1$-fibres, that is $t=[t_{1,0},\ldots,t_{1,r_1}]$ and
$d(t)=\sum_{i=0}^{r_1} (d(t_{1,i})+1) = r_1 + 1 + \sum_{i=0}^{r_1}
d(t_{1,i})$.

Let $F$ be an $\epin^X$-module. For fixed $1\leq j\leq n$ and $i\in[r_j]$
such that $f_j(i) = f_j(i+1)$ or
$j=1$  we define the map $d_i^j$
$$
F(t,\phi) \stackrel{d_i^j}{\lra}
 \bigoplus_{(t',\phi') =(X \stackrel{\phi'}{\rightarrow}[r_n]
\stackrel{g_n}{\rightarrow}\ldots \stackrel{g_{j+2}}{\rightarrow}
[r_{j+1}] \stackrel{d_if_{j+1}}{\rightarrow} [r_j-1]
\stackrel{f_j|_{i=i+1}}{\rightarrow} \ldots
\stackrel{f_2}{\rightarrow} [r_1]) \in \epin^X} F(t',\phi')$$ as
\begin{equation} \label{eq:dijX}
\begin{cases}
d_i^j = \sum_{\tau_{j+1}^{i,j} \in
\mathrm{Sh}(f_{j+1}^{-1}(i),f_{j+1}^{-1}(i+1))}
\varepsilon(\tau_{j+1}^{i,j};t_{j,i},t_{j,i+1})
F(\tau_n^{i,j},\ldots, \tau_{j+1}^{i,j}, d_i,\id, \ldots, \id)&
\text{ if }  j<n, \\ 
d_i^n=\epsilon(t_{n,i};t_{n,i+1})F( d_i,\id, \ldots, \id)& \text{ if
} j=n.
\end{cases}
\end{equation}

\begin{itemize}
\item[]
\item
If $F$ is an $\epin^X$-module, then the \emph{$(E_n,X)$-chain
complex of $F$} is the $n$-fold chain complex whose
$(r_n,\ldots,r_1)$ spot is
\begin{equation}
  \label{eq:cenX}
  C^{E_n,X}_{(r_n,\ldots,r_1)}(F) =
\bigoplus_{X\stackrel{\phi}{\rightarrow}[r_n]
\stackrel{f_n}{\rightarrow} \ldots \stackrel{f_2}{\rightarrow} [r_1]
\in \epin^X} F(X\stackrel{\phi}{\rightarrow}[r_n]
\stackrel{f_n}{\rightarrow}  \ldots \stackrel{f_2}{\rightarrow}
[r_1]).
\end{equation}
The differential in the $j$-th coordinate is
$$ \partial_j \colon C^{E_n,X}_{(r_n,\ldots,r_j, \ldots, r_1)}(F)
\rightarrow C^{E_n,X}_{(r_n,\ldots,r_j-1, \ldots, r_1)}(F)$$ with

$$ \partial_j := \sum_{i|f_j(i) = f_j(i+1)} (-1)^{s_{j,i}} d_i^j,$$
 where the $s_{j,i}$ are  obtained as follows: drawing the tree $t$ on a
plane with its
 root at the bottom, one can label its edges from bottom to top and left
to right; the integer $s_{j,i}$ is the sum of the label of the right
most top edge of the tree $t_{j,i}$ and the degrees of the elements
in $X$ which are in the fibre of the leaves that are to the left of
the top edge so defined, including it.

\item
The \emph{$(E_n,X)$-homology of the $\epin^X$-module $F$}, $H_*^{E_n,X}(F)$ is
defined to be the homology of the total complex associated to \eqref{eq:cenX}.
\end{itemize}
\end{defn}
\begin{lem} \label{lem:chain2}
The $k$-modules $C^{E_n,X}_{(r_n,\ldots,r_1)}(F)$ constitute an
$n$-fold chain complex.
\end{lem}
The proof of the lemma follows after example \ref{ex:tree14}.
\begin{ex}\label{ex:tree14}

Let $t$ be the following tree of degree 14 with its edges labelled,
$X=\{a_1,\ldots,a_{13}\}$ and $\phi$ can be read off the picture.

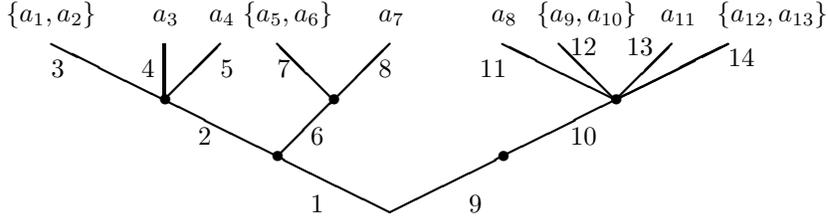
\raisebox{-3cm}{
\begin{picture}(10,10)
\setlength{\unitlength}{1.5cm}
\thicklines
\put(5,0){\line(2,1){3}}
\put(5,0){\line(-2,1){3}}
\put(4,0.5){\line(1,1){1}}
\put(4.5,1){\line(-1,1){0.5}}
\put(3,1){\line(1,1){0.5}}
\put(3,1){\line(0,1){0.5}}
\put(7,1){\line(-2,1){1}}
\put(7,1){\line(2,1){1}}
\put(7,1){\line(-1,1){0.5}}
\put(7,1){\line(1,1){0.5}}
\put(1.6,1.7){$\{a_1,a_2\}$}
\put(2.9,1.7){$a_3$}
\put(3.4,1.7){$a_4$}
\put(3.7,1.7){$\{a_5,a_6\}$}
\put(4.9,1.7){$a_7$}
\put(5.9,1.7){$a_8$}
\put(6.3,1.7){$\{a_9,a_{10}\}$}
\put(7.4,1.7){$a_{11}$}
\put(7.9,1.7){$\{a_{12},a_{13}\}$}
\put(4.3,0){$1$}
\put(3.3,0.6){$2$}
\put(3.95,0.45){$\bullet$}
\put(2.0,1.2){$3$}
\put(2.8,1.2){$4$}
\put(3.5,1.2){$5$}
\put(2.95,0.95){$\bullet$}
\put(4.3,0.6){$6$}
\put(4.0,1.2){$7$}
\put(4.9,1.2){$8$}
\put(4.45,0.95){$\bullet$}
\put(5.7,0){$9$}
\put(5.95,0.45){$\bullet$}
\put(6.6,0.6){$10$}
\put(6.95,0.95){$\bullet$}
\put(5.8,1.2){$11$}
\put(6.6,1.4){$12$}
\put(7.1,1.4){$13$}
\put(8,1.3){$14$}
\end{picture}
}
\medskip

This tree represents the object $X\stackrel{\phi}{\rightarrow}[8]
\stackrel{f_3}{\longrightarrow} [2] \stackrel{f_2}{\longrightarrow} [1] \in
\epi_3^X$, where the map $f_3$ maps $0,1,2$ to $0$, it sends $3,4$ to $1$ and
$5,6,7,8$ to $2$ and $f_2$ is $d_0$.

With our notation the tree $t_{1,0}$ is the $2$-level tree whose root is
the vertex above the edge labelled by 1, the tree $t_{1,1}$ is the subtree
above the edge with label 9, the tree $t_{2,0}$ is the $1$-level tree
above the label 2, $t_{2,1}$ the one above the label 6 and $t_{2,2}$ the
one above the label 10.

We have to determine the differentials $\partial_1,\partial_2$ and
$\partial_3$.

\medskip

In our example the differential $\partial_1$ glues the edges
labelled by $1$ and $9$ and shuffles the subtrees
$t_{1,0}=[t_{2,0},t_{2,1}]$ and $t_{1,1}= [t_{2,2}]$.  One has
$\partial_1=(-1)^{8+d(\{a_1,\ldots,a_7\})}d^1_0$ where $8$ is the
label of the right most edge of $t_{1,0}$. In addition we have the
shuffle signs. One has $d(t_{2,0})=3+d(\{a_1,\ldots,a_4\}),
d(t_{2,1})=2+d(\{a_5,a_6,a_7\})$,
$d(t_{1,0})=7+d(\{a_1,\ldots,a_7\})$ and
$d(t_{2,2})=4+d(\{a_8,\ldots,a_{13}\})$. In the expansion of $d^1_0$
there are 3 shuffles involved: $\id, (132)$, and $(231)$, written in
image notation. The first is coming with sign $+1$, the second one
with sign $(-1)^{(d(t_{2,1})+1)(d(t_{2,2})+1)}$ and the third one
with sign $(-1)^{(d(t_{2,0})+1+d(t_{2,1})+1)(d(t_{2,2})+1)}$. For
instance the image of the latter shuffle is in $F((t',\phi'))$ where
$(t',\phi')$ is the following tree:

\medskip

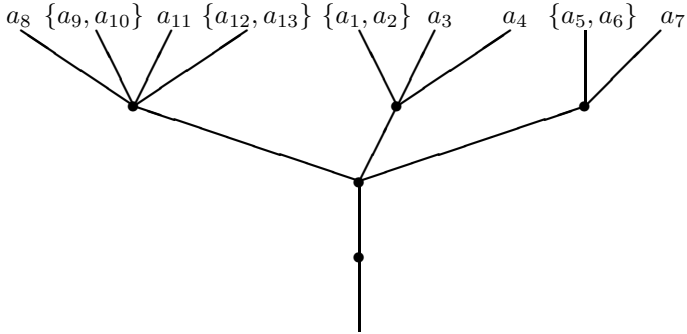
\raisebox{-5cm}{
\begin{picture}(10,10)
\setlength{\unitlength}{1cm}
\thicklines
\put(5,0){\line(0,1){2}}
\put(5,2){\line(-3,1){3}}
\put(5,2){\line(1,2){0.5}}
\put(5,2){\line(3,1){3}}
\put(2,3){\line(-3,2){1.5}}
\put(2,3){\line(-1,2){0.5}}
\put(2,3){\line(1,2){0.5}}
\put(2,3){\line(3,2){1.5}}
\put(5.5,3){\line(-1,2){0.5}}
\put(5.5,3){\line(1,2){0.5}}
\put(5.5,3){\line(3,2){1.5}}
\put(8,3){\line(0,1){1}}
\put(8,3){\line(1,1){1}}
\put(4.9,0.9){$\bullet$}
\put(4.9,1.9){$\bullet$}
\put(1.9,2.9){$\bullet$}
\put(5.4,2.9){$\bullet$}
\put(7.9,2.9){$\bullet$}
\put(0.3,4.1){$a_8$}
\put(0.8,4.1){$\{a_9,a_{10}\}$}
\put(2.3,4.1){$a_{11}$}
\put(2.9,4.1){$\{a_{12},a_{13}\}$}
\put(4.5,4.1){$\{a_{1},a_{2}\}$}
\put(5.9,4.1){$a_3$}
\put(6.9,4.1){$a_4$}
\put(7.5,4.1){$\{a_{5},a_{6}\}$}
\put(9,4.1){$a_7$}
\end{picture}
}
\medskip

The differential $\partial_2$ is $(-1)^{5+d(\{a_1,\ldots,a_4\})}
d^2_0$ where $5$ is the label of the  right most top edge of
$t_{2,0}$. The shuffles involved in the computation  of $d^2_0$ are
the  $(3,2)$-shuffles. For such a $(3,2)$-shuffle $\tau$ the
associated sign is given by $\epsilon(\tau;t_{2,0},t_{2,1})$ where
$t_{2,0}=[t_{3,0},t_{3,1},t_{3,2}]$ and $t_{2,1}=[t_{3,3},t_{3,4}]$.

The differential  $\partial_3$ is given by
$\partial_3=(-1)^{3+d(\{a_1,a_2\})}
d^3_0+(-1)^{4+d(\{a_1+a_2+a_3\})} d^3_1 +
(-1)^{7+d(\{a_1,\ldots,a_6\})} d^3_3+
(-1)^{11+d(\{a_1,\ldots,a_8\})}
d^3_5+(-1)^{12+d(\{a_1,\ldots,a_{10}\})} d^3_6+(-1)^{13+
d(\{a_1,\ldots,a_{11}\})}d^3_7.$
\end{ex}

\begin{proof}[Proof of \ref{lem:chain1} and \ref{lem:chain2}]
The proof that $d=\sum_j\partial_j$ satisfies $d^2=0$ is done by
induction on  $n$. Since the expression of $d$ in \ref{def:enchains}
coincides with the one in  \ref{def:enchainsX} when $d(X)=0$ it is
enough to prove lemma \ref{lem:chain2}.

Let  $(t,\phi):X \stackrel{\phi}{\rightarrow} [r_n]
\stackrel{f_n}{\rightarrow} \ldots \stackrel{f_{j+1}}{\rightarrow}
[r_j] \stackrel{f_j}{\rightarrow} \ldots \stackrel{f_2}{\rightarrow}
[r_1]$ be an $(X,n)$-level tree in $\epi^X_n$. For $j<n$ and  for
$i$ such that $f_j(i)=f_j(i+1)$ (we denote by $S^t_j$ the set of
such $i$'s) and for
$\tau\in\Sh(f_{j+1}^{-1}(i),f_{j+1}^{-1}(i+1))=\Sh^t_{i,j}$ we denote
by $\tau_{i,j}$ the map in $\epin^X$ defined in lemma \ref{lemma:di}
\begin{equation*}
\xymatrix{X\ar[r]^{\phi}\ar[d]_{id}&
{[r_n]} \ar[r]^{f_n} \ar[d]_{\tau_{n}^{i,j}}& {[r_{n-1}]} \ar[r]^{f_{n-1}}
  \ar[d]_{\tau_{n-1}^{i,j}}&
  {\ldots}\ar[r]^{f_{j+2}}&[r_{j+1}] \ar[d]_{\tau_{j+1}^{i,j}}
  \ar[r]^{f_{j+1}} & [r_j]  \ar[r]^{f_{j}} \ar[d]_{d_i} & 
{[r_{j-1}]} \ar[d]_{\id}
  \ar[r]^{f_{j-1}} & {\ldots} \ar[r]^{f_2}  &  {[r_1]} \ar[d]_{\id}\\
X\ar[r]^{\tau_{n}^{i,j}\phi}&{[r_n]} \ar[r]^{g_n^\tau} & {[r_{n-1}]}
\ar[r]^{g_{n-1}^\tau} & 
  {\ldots}\ar[r]^{g_{j+2}^\tau}& [r_{j+1}]\ar[r]^{d_if_{j+1}} &
  [r_j-1]  \ar[r]^{f_{j}|_{i=i+1}} & {[r_{j-1}]} 
\ar[r]^{f_{j-1}} &{\ldots}   \ar[r]^{f_2} &  {[r_1].}
}
\end{equation*}
For $j=n$ and $i$ such that $f_n(i)=f_n(i+1)$ there are no shuffles
involved  and we denote by
$d_{i,n}$ the corresponding map. We have to prove that
$\partial=\sum_j\partial_j$, where $\partial_j=\sum_{i|f_j(i)=f_j(i+1)}
(-1)^{s_{j,i}}d_i^j$ and $d_i^j$ is defined by relation (\ref{eq:dijX}),
satisfies $\partial^2=0$.
For $x\in F(t,\phi)$, $\partial(x)$ has the following form
\begin{equation}\label{eq:partialtree}
\partial(x)=\sum_{j<n,i\in S^t_j,\tau\in\Sh^t_{i,j}}\pm F(\tau_{i,j})(x)+
\sum_{i\in S^t_n}\pm F(d_{i,n})(x),
\end{equation}
where the signs are described in definition \ref{def:enchainsX}.

Let  $T_n^X$ be the free $k$-module generated by the trees in
$\epin^X$ and for a generator $t \in T_n^X$ consider the map
$$\partial(t)=\sum_{j<n,i\in S^t_j,\tau\in\Sh^t_{i,j}}\pm \tau_{i,j}(t)+
\sum_{i\in S^t_n}\pm d_{i,n}(t).$$
Then $\partial^2(t)=0$ in $T_n^X$ implies that $\partial^2(x) = 0$
for all $x \in F(t,\phi)$.

In order to prove that $\partial^2(t)=0$ for any $t\in\epin^X$ we
use the construction of the iterated bar construction given by
Eilenberg and Mac Lane in \cite[sections 7--9]{EM}. We differ from
their convention by using the left instead of the right action of
the symmetric group. If $(A,\partial)$ is a differential graded
commutative algebra then $BA$ is a differential graded commutative
algebra with a differential that is the sum of a residual boundary
$$\partial_r([a_1,\ldots,a_k])=
\sum_{i=1}^k(-1)^{i+d(a_1)+...+d(a_{i-1})}[a_1,\ldots,\partial
a_i,\ldots, a_k]$$ and a simplicial boundary
$$\partial_s([a_1,\ldots,a_k])=\sum_{i=1}^{k-1}(-1)^{i+d(a_1)+...+d(a_{i})}[a_1,
\ldots,a_i\cdot  a_{i+1},\ldots, a_k].$$ The graded commutative
product of $a=[a_1,\ldots,a_k]$ and $b=[a_{k+1},\ldots,a_{k+l}]$  is
given by the shuffle product
$$[a_1,\ldots,a_k]*[a_{k+1},\ldots,a_{k+l}]=
\sum_{\sigma\in\Sh(k,l)}\varepsilon(\sigma;a,b)
[a_{\sigma^{-1}(1)},\ldots,a_{\sigma^{-1}(k+l)}]$$ where
$\varepsilon(\sigma;a,b)$ picks up a factor
$(-1)^{(d(a_i)+1)(d(a_j)+1)}$  whenever
$\sigma(i)>\sigma(j)$ but $i<j$.

\medskip

Assume that $n=1$. To an ordered finite set $X$ of graded elements
one can associate a graded algebra $A = \bigoplus_{I\subset X} k[e_I]$
with $d(e_I)=d(I)$ and with the multiplication given by
$$e_Ie_J=\begin{cases} \epsilon(I;J) e_{I\sqcup J}, & \rm{if\ }I\cap
  J=\varnothing \\ 0, & \rm{if\ not}.\end{cases}$$
The algebra $A$ is graded commutative thanks to relation
(\ref{eq:epsilonAB}). Since the algebra $A$ has no differential the
boundary in $BA$ reduces to the simplicial boundary which is

$$\partial_s([e_{X_0},\ldots,e_{X_k}])=\sum_{i=0}^{k-1}(-1)^{i+1+d(X_0)+...+d(X_{i})}
[e_{X_0},\ldots,e_{X_i}\cdot  e_{X_{i+1}},\ldots, e_{X_k}].$$

If we assume furthermore that the family $(X_0,\ldots,X_k)$ forms a
partition of $X$ then $\partial_s([e_{X_0},\ldots,e_{X_k}])$ is a
sum of elements satisfying the same property. As a consequence the
free $k$-module generated by the trees in $\epi_1^X$ is a subcomplex
of $BA$ and $\partial_1=\partial_s$ is a differential, since
$s_{1,i}$ is precisely $i+1+\sum_{j=0}^i d(X_j)$.

\medskip

Assume $n>1$. We denote by $T_n^X$ the free $k$-module generated by
the $(X,n)$-level trees. We denote by $T_n$ the free $k$-module generated
by the $n$-level trees whose top vertices are labelled by subsets of
$X$. The labels  $X_a$ and $X_b$ of two different vertices $a$ and $b$
may have a non-empty intersection. 
An $(X,n)$-level
tree $t=[t_{1,0},\ldots,t_{1,r_1}]$ is considered as an element of
$BT_{n-1}$. The $k$-module $T_{n-1}$ is endowed with the shuffle
product since $T_{n-1}=BT_{n-2}$. More precisely, for
$a=[t_1,\ldots, t_p]$ and $b=[t_{p+1},\ldots,t_{p+q}]$ two
$(n-1)$-level trees where $t_i$ is an $(n-2)$-level tree for every
$i$, one has
$a*b=\sum_{\sigma\in\Sh_{p,q}}\epsilon(\sigma;a,b)[t_{\sigma^{-1}(1)},\ldots,
t_{\sigma^{-1}(p+q)}]$. 
From Eilenberg and Mac Lane \cite[7.1]{EM} one gets that
$(BT_{n-1},\partial_s+\partial_r)$ is a complex of $k$-modules.

Let $t$ be an $(X,n)$-level tree. One has
$t=[t_{1,0},\ldots,t_{1,r_1}]$ where each $t_{1,i}$ is labelled by
$X_i\subset X$ and the family $(X_i)_{0\leq i\leq r_1}$ is a
partition of $X$. Hence $T_n^X$ is a $k$-submodule of $BT_{n-1}$.
The simplicial boundary preserves $T_n^X$ for $t_{1,i}*t_{1,i+1}\in
T_{n-1}^{X_i\cup X_{i+1}}$. The residual boundary preserves $T_n^X$
by induction on $n$. Furthermore
\begin{multline*}
\partial_s(t)=\sum_{i=0}^{r_1-1} (-1)^{i+1+d(t_{1,0})+\ldots
  d(t_{1,i})}[t_{1,0}, \ldots,t_{1,i}*t_{1,i+1},\ldots,
t_{1,r_1}]=\\
\sum_{i=0}^{r_1-1}
(-1)^{s_{1,i}}\sum_{\tau\in\Sh_{i,1}^t}\epsilon(\tau;t_{1,i},
t_{1,i+1})\tau_{i,1}(t)=\partial_1(t). 
\end{multline*}
Similarly, by induction on $n$ we prove that the sum
$\partial_2+\ldots +\partial_n$ corresponds to the residual
boundary. As a consequence $T_n^X$ is a sub-complex of $BT_{n-1}$
and one gets that $\partial^2=0$.

\end{proof}

As an example, we will determine the zeroth $E_n$-homology of an
$\epin$-functor  $F$. In total degree zero there is just one
summand, namely $F([0]\stackrel{\id_{[0]}}{\lra} \ldots
\stackrel{\id_{[0]}}{\lra} [0])$. The modules
$C^{E_n}_{(0,1,0,\ldots,0)}(F)$, $\ldots$,
$C^{E_n}_{(0,\ldots,0,1)}(F)$ are all trivial, so the only boundary
term that can occur is caused by the unique map
$$ C^{E_n}_{(1,0,\ldots,0)}(F) \lra C^{E_n}_{(0,\ldots,0)}(F).$$
Therefore
\begin{equation}  \label{eq:H0}
 H^{E_n}_0(F) \cong F([0]\stackrel{\id_{[0]}}{\lra} \ldots
\stackrel{\id_{[0]}}{\lra}
[0])/\mathrm{image}(F([1]\stackrel{d_0}{\lra} [0]
\stackrel{\id_{[0]}}{\lra} \ldots \stackrel{\id_{[0]}}{\lra} [0])).
\end{equation}

We can view an $\epin$-module $F$ as an $\mathrm{Epi}_k$-module for
all $k \leq n$ via the functors $\iota_n^k$.
\begin{prop} \label{prop:epikepinhomology}
For every $\epin$-module $F$ there is a map of chain complexes $
\mathrm{Tot}(C_*^{E_k}(F \circ \iota_n^k))  \lra
\mathrm{Tot}(C_*^{E_n}(F))$ and therefore a map of graded
$k$-modules
$$ H_*^{E_k}(F  \circ \iota_n^k) \lra H_*^{E_n}(F). $$
\end{prop}
\begin{proof}
There is a natural identification of the module
$C_{(r_k,\ldots,r_1)}^{E_k}(F \circ \iota_n^k)$ with the module
$C_{(r_k,\ldots,r_1,0,\ldots,0)}^{E_n}(F)$ and this includes
$\mathrm{Tot}(C_*^{E_k}(F \circ \iota_n^k))$ as a subcomplex into
$\mathrm{Tot}(C_*^{E_n}(F))$.
\end{proof}

\subsection{Relationship to higher order Hochschild homology}
\label{subsec:comp} 
For a non-unital commutative $k$-algebra $\bar{A}$
we define $\L^n(\bar{A})\colon \epin \rightarrow \kmod$ as
$$ \L^n(\bar{A})([r_n] \stackrel{f_n}{\rightarrow}  \ldots
\stackrel{f_2}{\rightarrow} [r_1]) = \bar{A}^{\otimes (r_n+1)}.$$
A morphism
$$ \xymatrix{{[r_n]} \ar[d]^{\sigma_n} \ar[r]^{f_n} & {[r_{n-1}]}
  \ar[d]^{\sigma_{n-1}} \ar[r]^{f_{n-1}} &
  {\ldots}  \ar[r]^{f_2} & {[r_1]} \ar[d]^{\sigma_{1}} \\
{[r'_n]} \ar[r]^{f'_n} & {[r'_{n-1}]} \ar[r]^{f'_{n-1}} &
  {\ldots} \ar[r]^{f'_2} & {[r'_1]}}
$$
induces a map $\bar{A}^{\otimes (r_n+1)} \rightarrow \bar{A}^{\otimes
  (r'_n+1)}$ via
$$ a_0 \otimes \ldots \otimes a_{r_n} \mapsto (\sigma_n)_*(a_0 \otimes
\ldots \otimes a_{r_n}) = b_0 \otimes \ldots \otimes b_{r'_n}$$ with
$b_i = \prod_{\sigma_n(j) = i} a_j$. The $E_n$-homology of the
functor $\L^n(\bar A)$ coincides with the homology of the $n$-fold
 bar construction of the non-unital algebra $\bar A$, hence with the
$E_n$-homology of $\bar A$. The total complex has been described  in
\cite[Appendix]{F} and it coincides with ours.

There is a correspondence between augmented commutative $k$-algebras
and non-unital $k$-algebras that sends an augmented $k$-algebra $A$
to its augmentation ideal $\bar{A}$. Under this correspondence, the
$(m+n)$-th homology group of the $n$-fold bar construction $B^n(A)$
is isomorphic to the $m$-th homology group of the $n$-fold iterated
bar construction of $\bar{A}$,  $B^n(\bar{A})$. As the chain complex
$B(A)$ is the chain complex for the Hochschild homology of $A$ with
coefficients in $k$ (compare \cite[(7.5)]{EM}), we can express
$B(A)$ as $A \tilde{\otimes} \mathbb{S}^1$. Here, $\mathbb{S}^1$ is
the simplicial model of the $1$-sphere, which has $n+1$ elements in
simplicial degree $n$ and $(A \tilde{\otimes} \mathbb{S}^1)_n = k
\otimes A^{\otimes n}$. Therefore
$$B^n(A) \cong (\ldots (A \tilde{\otimes} \mathbb{S}^1) \ldots )
\tilde{\otimes} \mathbb{S}^1 \cong A \tilde{\otimes}
((\mathbb{S}^1)^{\wedge n}) \cong A \tilde{\otimes} \mathbb{S}^n$$
which gives rise to higher order Hochschild homology of order $n$ of
$A$ with coefficients in $k$, $HH^{[n]}_*(A;k)$,  in the sense of
Pirashvili \cite{P}. Thus, $HH_{*+n}^{[n]}(A;k) \cong
H_*^{E_n}(\bar{A})$.

By proposition \ref{prop:epikepinhomology} there is a sequence of
maps
\begin{equation} \label{eq:sequence}
HH_{*+1}(A;k) \cong H_*^{\mathrm{bar}}(\bar{A}) = H^{E_1}_*(\bar{A})
\ra H^{E_2}_*(\bar{A}) \ra H^{E_3}_*(\bar{A}) \ra \ldots
\end{equation}
and the map from $H^{E_1}_*(\bar{A})$ to the higher $E_n$-homology
groups is given on the chain level by the inclusion of
$C_m^{\mathrm{bar}}(\bar{A})$ into
$C^{E_n}_{(m,0,\ldots,0)}(\bar{A})$.

Suspension induces maps
$$
\xymatrix{
{HH_\ell(A;k) = \pi_\ell\L(A;k)(\mathbb{S}^1)} \ar[dr] \ar[r] &
{HH_{\ell +1}^{[2]}(A;k) =  \pi_{\ell+1}\L(A;k)(\mathbb{S}^2)} \ar[d] \ar[r]
& {\ldots} \\
& {H\Gamma_{\ell-1}(A;k) \cong \pi_\ell^s(\L(A;k)).}& }$$
For the last isomorphism see \cite{PR2}. Fresse proves a comparison
\cite[9.6]{F} between Gamma homology of $A$ and $E_\infty$-homology of
$\bar{A}$. Using the isomorphisms above
this sequence gives rise to a sequence of maps involving graded $k$-modules
that are isomorphic to the ones in \eqref{eq:sequence}.

The explicit form of the suspension maps is described in
\cite[(7.9)]{EM}: an element $a \in \bar{A}$ is sent to $[a]$ in the
bar construction. The iterations of this map correspond precisely to
the maps $\iota^{n-1}_n\colon B^{n-1}(\bar{A}) \rightarrow
B^n(\bar{A})$. Therefore we actually have an isomorphism of
sequences, \ie, the suspension maps $HH^{[n]}_{\ell+n}(A;k)
\rightarrow HH_{\ell +n+1}^{[n+1]}(A;k)$ are
  related to the natural maps $H^{E_n}_\ell(\bar{A}) \rightarrow
H^{E_{n+1}}_\ell(\bar{A})$ via the isomorphisms
$HH_{*+n}^{[n]}(A;k) \cong H_*^{E_n}(\bar{A})$.

We have $\L^n(\bar{A})([1]\stackrel{\id_{[0]}}{\lra} \ldots
\stackrel{\id_{[0]}}{\lra} [0]) = \bar{A}^{\otimes 2}$ and hence for
all $n \geq 1$ the zeroth $E_n$-homology group is
$$ H^{E_n}_0(\bar{A}) \cong \bar{A}/\bar{A}\cdot \bar{A}.$$

As a concrete example, we calculate the $E_n$-homology of a polynomial
algebra in one and in two variables.  Let $\Gamma$ denote the skeleton
of the category of finite pointed sets with objects $[n] =
\{0,\ldots,n\}$ for $n \geq 0$ with $0$ as basepoint. Pirashvili's
definition of
$n$-th order Hochschild homology is given for $\Gamma$-modules, \ie,
functors from $\Gamma$ to $k$-modules. In particular, we consider
$\mathcal{L}(A;k)\colon \Gamma \rightarrow k\mathrm{-mod}$, for an
augmented commutative algebra $A$ with $\mathcal{L}(A;k)[n] =
A^{\otimes n}$ where maps in $\Gamma$ induce multiplication in $A$,
augmentation to $k$ or insertion of units.
The functor $\mathcal{L}$ for a polynomial
algebra with coefficients in $k$ evaluated on a simplicial model of
the $n$-sphere is the symmetric algebra functor evaluated on the
$n$-sphere and thus arguing as in \cite[section 4]{RiRo} we obtain
$$ H^{E_n}_*(\overline{k[x]}) \cong HH_{*+n}^{[n]}(k[x];k) =
H_{*+n}(\mathcal{L}(k[x];k)(\mathbb{S}^n)) \cong
H_{*+n}(\mathrm{Sym} \circ L(\mathbb{S}^n)) \cong
H_{*+n}(SP(\mathbb{S}^n);k).$$

Here, $SP$ stands for the infinite symmetric product and $L$ is the
$\Gamma$-module that sends $[n]$ to the free $k$-module generated by
the set $\{1,\ldots,n\}$. Note that in this case $k[x]$ is augmented
over $k$ via the augmentation $\varepsilon_0$ that sends $x$ to
zero. The augmentation affects the $k[x]$-module structure of $k$,
but in \cite[4.1]{RiRo} it is shown that the resulting homotopy
groups are independent of the module structure.

Evaluated on an $n$-sphere, the functor $SP$ yields an
Eilenberg-MacLane space of type $(\mathbb{Z},n)$ and hence the above
is  isomorphic to $H_{*+n}(K(\mathbb{Z},n);k)$.

Consider the $\Gamma$-module $\L(k[x,y];k))$. To $[n]$ it associates $k
\otimes k[x,y]^{\otimes n} \cong k[x]^{\otimes n} \otimes
k[y]^{\otimes n} \cong \L(k[x];k))[n] \otimes \L(k[y];k))[n]$. A
morphism of finite pointed sets $f\colon [n] \rightarrow [m]$ sends
$\lambda \otimes a_1 \otimes \ldots \otimes a_n$ (with $\lambda \in
k$ and $a_i$ in $k[x,y]$) to $\mu \otimes b_1 \otimes \ldots \otimes
b_m$ where $b_i = \prod_{f(j) = i} a_j$ and $\mu = \lambda \cdot
\prod_{f(j)=0, j\neq 0} \varepsilon(a_j)$. Therefore the above 
isomorphism of $\L(k[x];k))[n] \otimes \L(k[y];k))[n]$  and
$\L(k[x,y];k))[n]$ induces an isomorphism of $\Gamma$-modules
between $\L(k[x,y];k))(\mathbb{S}^n)$ and the pointwise tensor
product $\L(k[x];k))(\mathbb{S}^n) \otimes
\L(k[y];k))(\mathbb{S}^n)$. Furthermore, we get
\begin{align*}
\pi_*((\mathrm{Sym} \circ L(\mathbb{S}^n)) \otimes (\mathrm{Sym} \circ
L(\mathbb{S}^n))) & \cong \pi_*(\mathrm{Sym} \circ (L(\mathbb{S}^n)
\oplus
L(\mathbb{S}^n)) \\
& \cong \pi_*(\mathrm{Sym} \circ (L(\mathbb{S}^n \vee \mathbb{S}^n))
\cong H_*(SP(\mathbb{S}^n \vee \mathbb{S}^n);k) \\
& \cong H_*(K(\mathbb{Z} \times \mathbb{Z},n);k).
\end{align*}

Therefore, we obtain that
$$ H^{E_n}_*(\overline{k[x,y]}) \cong H_{*+n}^{[n]}(k[x,y]) \cong
H_{*+n}(K(\mathbb{Z} \times \mathbb{Z},n);k).$$

We close this part with the description of a spectral sequence.
Our description of $E_2$-homology leads to the following result.
\begin{prop}
If $\bar{A}$ and $H_*^{\mathrm{bar}}(\bar{A})$ are $k$-flat, then
there is a spectral sequence
$$ E^1_{p,q} = \bigoplus_{\ell_0 + \ldots + \ell_q =p- q}
H_{\ell_0}^{\mathrm{bar}}(\bar{A}) \otimes \ldots \otimes
H_{\ell_q}^{\mathrm{bar}}(\bar{A}) \Rightarrow
H^{E_2}_{p+q}(\bar{A}) $$ where the $d_1$-differential is induced by
the shuffle differential.
\end{prop}
\begin{proof}
The double complex for $E_2$-homology looks as follows:
$$ \xymatrix{
{} & {} & {\vdots} \ar[d] & {}\\
{} & {} & {\bar{A}^{\otimes 3}} \ar[d] & {\ldots}\ar[l]\\
{} & {\bar{A}^{\otimes 2}} \ar[d] & {\bar{A}^{\otimes 3}} \ar[l]
\ar[d] &  {\ldots}\ar[l]\\
{\bar{A}} & {\bar{A}^{\otimes 2}} \ar[l] & {\bar{A}^{\otimes 3}}
\ar[l] &  {\ldots}\ar[l] }$$ The horizontal maps are induced by the
$b'$-differential whereas the vertical maps are induced by the
shuffle maps. The horizontal homology of the bottom row is precisely
$H_*^{\mathrm{bar}}(\bar{A})$. We can interpret the second row as
the total complex associated to the following  double complex:
$$ \xymatrix{
{\vdots} \ar[d]_{\id \otimes b'} & {\vdots} \ar[d]_{\id \otimes b'} &
{\vdots} \ar[d]_{\id \otimes b'}  & {}\\
{\bar{A} \otimes \bar{A}^{\otimes 3}} \ar[d]_{\id \otimes b'} &
{\bar{A}^{\otimes 2} 
  \otimes \bar{A}^{\otimes 3}} \ar[d]_{\id \otimes b'}
\ar[l]_{b'\otimes \id} &  {\bar{A}^{\otimes 3}
  \otimes \bar{A}^{\otimes 3}} \ar[d]_{\id \otimes b'} \ar[l]_{b'\otimes \id} &
{\ldots}\ar[l]_(0.35){b'\otimes \id} \\
{ \bar{A} \otimes \bar{A}^{\otimes 2}} \ar[d]_{\id \otimes b'} &
{\bar{A}^{\otimes 2}\otimes
  \bar{A}^{\otimes 2} } \ar[d]_{\id \otimes b'} \ar[l]_{b'\otimes \id}
&  {\bar{A}^{\otimes 3} \otimes
  \bar{A}^{\otimes 2}}  \ar[l]_{b'\otimes \id} \ar[d]_{\id \otimes b'}
&   {\ldots}\ar[l]_(0.35){b'\otimes \id}\\
{\bar{A} \otimes \bar{A}} & {\bar{A}^{\otimes 2} \otimes \bar{A} }
\ar[l]_{b'\otimes \id} &  {\bar{A}^{\otimes 3} \otimes\bar{A}}
\ar[l]_{b'\otimes \id} &  {\ldots}\ar[l]_(0.35){b'\otimes \id} }$$
Therefore the horizontal homology groups of the second row are the
homology of the tensor product of the
$C^{\mathrm{bar}}(\bar{A})$-complex with itself. Our flatness
assumptions guarantee that we obtain
$H_*^{\mathrm{bar}}(\bar{A})^{\otimes 2}$ as homology. An induction
then finishes the proof.
\end{proof}

\section{Tor interpretation of $E_n$-homology} \label{sec:comp}

The purpose of this section is to prove the following theorem

\begin{thm}\label{thm:mainthm}
For any $\epin$-module $F$
\begin{equation*}
H_p^{E_n}(F) \cong \mathrm{Tor}^{\epin}_p(\bepi_n, F), \text{ for all } p \geq
0
\end{equation*}
where
$$ \bepi_{n}(t) \cong \left\{
  \begin{array}{cr}
    k, & \text{ for } t=[0]\stackrel{\id_{[0]}}{\lra} \ldots
\stackrel{\id_{[0]}}{\lra} [0], \\
 0, & \text{ for } t \not=[0]\stackrel{\id_{[0]}}{\lra} \ldots
\stackrel{\id_{[0]}}{\lra} [0].
  \end{array} \right. $$
\end{thm}

As  will be explained in the proof below it is enough to prove that
for every  $t\in\epin$,
the representable functor $\epin^t=k[\epin(t,-)]$ is acyclic with respect
to  $E_n$-homology. In order to prove this we proceed by induction
on $n$.  The case $n=1$ has been treated in section
\ref{sec:barhomology}. We consider the  case
$n=2$  from proposition \ref{prop:Hd2} to corollary
\ref{cor:fork2}  and this case is the core of the proof. The induction
process
is  explained in proposition \ref{prop:Hdn}.
Unfortunately, the signs involved in the differential of the complex
$C_*^{E_n}(\epin^t)$ are not compatible with an induction process as
in  proposition
\ref{prop:Hdn} and this is the reason why we introduced the category
$\epin^X$  and the homology $H_*^{E_n,X}$ of an $\epin^X$-module
in definition \ref{def:enchainsX}. We explain in the proof below how
we use information on ${\epin^X}^{t,\phi}$ in order to get information on
$\epin^t$. Here $X$ is a finite
ordered set of  graded elements and ${\epin^X}^{t,\phi}$ is the
representable  functor $k[\epin^X[(t,\phi),-)]$.
Consequently propositions and lemmas \ref{prop:Hd2} to
\ref{cor:fork2}  are expressed in terms of the representable functors
${\epin^X}^{t,\phi}$ as well as the representable functors $\epin^t$,
whereas  proposition \ref{prop:Hdn} gives results on
${\epin^X}^{t,\phi}$ only.

\begin{proof}
As in the proof of proposition \ref{prop:n=1},
we have to show that $H_*^{E_n}(-)$ maps short exact sequences of
$\epin$-modules to long exact sequences, that $H_*^{E_n}(-)$
vanishes on projectives in  positive degrees and that $H_0^{E_n}(F)$ and
$\bepi_n \otimes_{\epin} F$ agree for all $\epin$-modules $F$. The
homology  $H_*^{E_n}(-)$ is the homology of a total complex
$C_*^{E_n}(-)$  sending short exact sequences as in \eqref{eq:ses} to
short  exact
sequences of chain complexes and therefore the first claim is true.
Note that the left $\epin$-module $\bepi_n$ is the cokernel of the map
between  contravariant representables
$$ (d_{0})_{*}\colon {\epin}_{,[1] \lra [0] \lra \ldots
\lra [0]} \rightarrow {\epin}_{,[0] \lra [0] \ldots \lra [0]}.$$
This remark together with the computation of $H_0^{E_n}(F)$ in
relation (\ref{eq:H0})  implies the last  claim, similarly to the
proof of  proposition \ref{prop:n=1}.

In order to show that $H_*^{E_n}(P)$ is trivial in positive degrees
for any projective $\epin$-module $P$ it suffices to show that the
representables $\epin^{t}$ are acyclic for any planar tree $t=[r_n]
\stackrel{f_n}{\rightarrow} \ldots \stackrel{f_2}{\rightarrow}
[r_1]$.

\medskip

Let $t$ be such an $n$-level tree, let $X$ be a finite ordered set and
let $\phi:X\rightarrow [r_n]$ be a
fixed surjection. Assume that every element in $X$ has degree
$0$. Then we claim that the
complexes
$C^{E_n}_*(\epin^t)$ and $C^{E_n,X}_*({\epin^X}^{t,\phi})$ are isomorphic. One has
$$k[\epin(t,t')] \cong \bigoplus_{\phi'\colon X \rightarrow
  [r'_n]}k[\epin^X((t,\phi),(t',\phi'))],$$
because any morphism of $n$-trees $\sigma\colon t \rightarrow t'$ determines a
component  $\phi' = \sigma_n \circ \phi$. This
defines an injective map $ k[\epin(t,t')]\rightarrow
\bigoplus_{\phi'\colon X \rightarrow
  [r'_n]}k[\epin^X((t,\phi),(t',\phi'))]$. As every morphism from $(t,\phi)$ to
$(t',\phi')$ is a morphism of  $n$-trees $\sigma\colon t
\rightarrow t'$ with $\sigma_n \circ \phi = \phi'$, the map is surjective.
By relations \eqref{eq:cen} and  \eqref{eq:cenX} one has

\begin{equation}\label{eq:Xnongraded}
C^{E_n}_*(\epin^t) = \bigoplus_{t'\in\epin}\epin(t,t') =
\bigoplus_{(t',\phi')\in\epin^X}
\epin^X((t,\phi),(t',\phi'))= C^{E_n,X}_*({\epin^X}^{t,\phi})
\end{equation}
and as every element of $X$ has degree zero, the differentials
$\partial_j$ coincide for all $j$.

In proposition \ref{prop:Hdn}, we will prove by induction that for an
$(X,n)$-level tree 
$(t,\phi): X \stackrel{\phi}{\rightarrow} [r_n]
\stackrel{f_n}{\rightarrow} \ldots \stackrel{f_{j+1}}{\rightarrow}
[r_j] \stackrel{f_j}{\rightarrow} \ldots \stackrel{f_2}{\rightarrow}
[r_1]$ with $X=\{x_0<\ldots<x_{r_n}\}$ and $\phi(x_i)=i$,
the representable ${\epin^X}^{t,\phi}$ is acyclic. In
particular, if every
element in $X$ has degree zero, then this implies that
$\epin^{t}$ is acyclic for any $n$-level tree $t$.

The case $n=1$ has been proved in proposition \ref{prop:n=1} in the
ungraded case and in remark \ref{R:gradedn=1} in the graded case.
For $n=2$ we study the bicomplex
$C^{E_2,X}_{(*,*)}({\epi_2^X}^{t,\phi})$. In proposition
\ref{prop:Hd2} we give  the $k$-module structure of the homology
with respect to the  differential $\partial_2$ and give its
generators in propositions \ref{prop:tophomology} and
\ref{prop:shufflecycles}. Corollaries  \ref{cor:nonfork2} and
\ref{cor:fork2}  state the result for $n=2$.  For the general case,
one uses  induction on $n$ and proposition \ref{prop:Hdn}. As a
consequence $H^{E_n}_*(\epi_n^t)=0$ for all $* \geq 0$ if
$t\not=[0]\lra[0]\ldots\lra [0]$ and in that case
$$H^{E_n}_*(\epi_n^{[0]\lra[0]\ldots\lra[0]})=\begin{cases} 0, &\mathrm{for}\
  *>0 \\ k, & \mathrm{for}\ *=0.\end{cases}$$
\end{proof}

In the following we need some technical tools from the homology of
small categories, as in Mitchell \cite[section 17]{M72}. We review
the standard resolution of a small category in the graded context.
References on the more general context of differential graded
categories can be found also in \cite{keller06}.

\begin{defn} \label{def:dgcat} Let $\mathcal C$ be a small category that
is a graded $k$-linear category, \ie, for every pair of objects
$a,b\in \mathcal C$ the set of morphisms $\mathcal C(a,b)$ is a
graded $k$-module and the structure maps of the category
$\mathcal{C}(b,c)\otimes\mathcal{C}(a,b)\rightarrow
\mathcal{C}(a,c)$ are morphisms of graded $k$-modules. We denote by
$\mathcal{C}\otimes \mathcal{D}$ the tensor product of two small
graded $k$-categories, which is given by the product on objects and
the tensor product on morphisms. The standard resolution of
$\mathcal C$ is the simplicial bifunctor from $\mathcal
C^{op}\otimes \mathcal C$ to the category of graded $k$-modules
defined by $$S_n(\mathcal C)=\bigoplus_{p_1,\ldots,p_{n+1}} \mathcal
C(-,p_1)\otimes\mathcal C(p_1,p_2)\otimes\ldots\otimes\mathcal
C(p_n,p_{n+1})\otimes\mathcal C(p_{n+1},-)$$ where the face maps are
$$d_i(\alpha_0\otimes\ldots\otimes\alpha_{n+1})
=(-1)^{d(\alpha_0)+\ldots+d(\alpha_i)} 
 \alpha_0\otimes\ldots\otimes\alpha_{i+1}\alpha_i \otimes \ldots
 \otimes \alpha_{n+1}$$ 
and the degeneracy maps are
$$s_j(\alpha_0\otimes\ldots\otimes\alpha_{n+1})=(-1)^{d(\alpha_0) + \ldots 
+d(\alpha_j)} \alpha_0\otimes\ldots\alpha_j\otimes\id\otimes\alpha_{j+1}\otimes
\ldots\otimes \alpha_{n+1}.$$
A covariant (contravariant) functor from the category
$\mathcal C$ to the category of graded $k$-modules is called a left
(right) graded $\mathcal C$-module. For any left graded $\mathcal
C$-module $L$ and any right graded $\mathcal C$ module $R$, we define

$$
H_*(R;L)=H_*(L\otimes_{\mathcal{C}^\mathrm{op}} S_*(\mathcal{C})
\otimes_{\mathcal{C}^\mathrm{op}} R)=\mathrm{Tor}^{\mathcal{C}}_*(R;L).$$
\end{defn}

Note that one can also use the normalized standard resolution,
combined with the Yoneda lemma to compute this homology
\begin{equation}\label{eq:normalized}
 L\otimes_{\mathcal{C}^\mathrm{op}} N_n(\mathcal{C})
\otimes_{\mathcal{C}^\mathrm{op}} R=\bigoplus_{p_1,\ldots,p_{n+1}}
L(p_1)\otimes\tilde{\mathcal C}(p_1,p_2)\otimes\ldots\otimes\tilde{\mathcal
C}(p_n,p_{n+1})\otimes R(p_{n+1}),
\end{equation}
where $\tilde{\mathcal C}(p_1,p_2)=\mathcal C(p_1,p_2)$ if
$p_1\not=p_2$ and $\tilde{\mathcal C}(p_1,p_2)$ is the cokernel of
the map $k\rightarrow\mathcal C(p_1,p_1)$ if $p_1=p_2$.

\smallskip

Let $a>0$ be an integer and let $Y$ be a graded ordered set and
$\pi=\pi^1\sqcup\ldots\sqcup \pi^a$ be a partition of $Y$. We denote by
$[a]^\pi$ the following graded category: objects in $[a]^\pi$ are the
elements $i$ for $0\leq i\leq a$ and morphisms in $[a]^\pi$ are the graded
$k$-modules given by

$$[a]^\pi(i,j)=\begin{cases} k \text{ in degree }
d(\pi^{i+1}\cup\ldots\cup\pi^j) ,& \text{ if } i\leq j, \\ 0 & \text{ if }
i>j, \end{cases}$$
where by convention $d(\varnothing)=0$.

The composition of  $\alpha\in [a]^\pi(j,k)$ and
$\beta\in [a]^\pi(i,j),\; i\leq j\leq k$ is given by
$$\alpha\circ\beta=\epsilon(\pi^{i+1}\cup\ldots\cup \pi^j
;\pi^{j+1}\cup\ldots\cup \pi^k)\alpha\beta.$$
The category $[a]^\pi$ is a poset category with a minimal element $0$ and
a maximal element $a$. We denote by $L_0$ the left graded $[a]^\pi$-module
which assigns $k$ in degree $0$ to $0$ and $0$ to $0<i\leq a$ and by $R_a$
the right graded $[a]^\pi$-module which assigns $k$ in degree $0$ to $a$
and $0$ to $0\leq i<a$.

\label{}
\begin{lem}\label{lem:torcomputation} For a tensor product
$[a_0]^{\pi_0}\otimes\ldots\otimes [a_{r_1}]^{\pi_{r_1}}$ of categories we
obtain
\begin{equation*} \mathrm{Tor}_n^{[a_0]^{\pi_0}\otimes\ldots\otimes
[a_{r_1}]^{\pi_{r_1}}}(R_{a_0}\otimes\ldots\otimes
R_{a_{r_1}};L_0\otimes\ldots\otimes L_0) \cong
\begin{cases} k, & \text{ if }
n=r_1+1 \text{ and } \forall j, a_j=1, \\
0, & \text{otherwise.}
\end{cases}
\end{equation*}
\end{lem}

\begin{proof}
The K\"unneth formula (see e.g. \cite[section 3.6]{We}) gives
\begin{multline*}
\mathrm{Tor}_n^{[a_0]^{\pi_0}\otimes\ldots\otimes
  [a_{r_1}]^{\pi_{r_1}}}(R_{a_0}\otimes\ldots\otimes 
R_{a_{r_1}};L_0\otimes\ldots\otimes L_0) \cong \\
\bigoplus_{n_0+\ldots+n_{r_1}=n}
\mathrm{Tor}_{n_0}^{[a_0]^{\pi_0}}(R_{a_0};L_0)\otimes\ldots\otimes
\mathrm{Tor}_{n_{r_1}}^{[a_{r_1}]^{\pi_{r_1}}}(R_{a_{r_1}};L_0).
\end{multline*}
Consequently it is enough to compute
$\mathrm{Tor}_{n}^{[a]^{\pi}}(R_{a};L_0)$. Relation
(\ref{eq:normalized}) gives the complex computing this homology: 
$$C_n([a]^\pi)=L_0\otimes_{\mathcal{C}^\mathrm{op}} N_n(\mathcal{C})
\otimes_{\mathcal{C}^\mathrm{op}} R_a=\begin{cases}  0, & \text { if } n=0, \\
k, & \text{ if } n=1, \\
\bigoplus_{0<p_2<\cdots<p_n<a} [a]^\pi(0,p_2)\otimes\ldots\otimes
[a]^\pi(p_n,a),  & \text{ if } n>1,\end{cases}$$
with the differential given by $d=\sum_{i=1}^{n-1} (-1)^i d_i$.

If $a=1$, then $C_n([a]^\pi)=0$ for $n>1$ and the result follows.

Assume $a>1$, and let $\xi_{i,j}$ be the generator
of $[a]^\pi(i,j)$ in degree $d(\pi^{i+1}\cup\ldots\cup\pi^j)$. The homotopy
$$h(\xi_{0p_2}\otimes\ldots\otimes\xi_{p_n a})=\begin{cases} 0, &
  \text{ if }  p_2=1, \\
(-1)^{d(\pi^1)} \varepsilon(\pi^1;\pi^2\cup\ldots\cup\pi^{p_2})
\xi_{01}\otimes \xi_{1p_2}\otimes\ldots\otimes\xi_{p_n a},  &\text{
if }  p_2>1,\end{cases}$$ proves that the complex is acyclic.
\end{proof}

\medskip

The $E_n$-homology of an $\epin$-module $F$ (resp.~the
$(E_n,X)$-homology of an $\epin^X$-module $F$) can be computed in
different ways,  since it is the homology of the total complex
associated to an $n$-complex. The notation  $H_*(F,\partial_i)$
stands for the homology of the complex $C_*^{E_n}(F)$ (resp.
$C_*^{E_n,X}(F)$) with respect to the differential  $\partial_i$.
The complex $(C_*^{E_n}(F),\partial_i)$ splits into subcomplexes

\begin{equation}\label{eq:splitcomplex}
C_{(s_n,s_{n-1},\ldots,s_{i+1},*,s_{i-1},\ldots,s_1)}^{E_n}(F)=
\bigoplus_{t=[s_n]\stackrel{g_n}{\lra}\ldots
\stackrel{g_{i+2}}{\lra}[s_{i+1} ]
\stackrel{g_{i+1}}{\lra}[*]\stackrel{g_{i}}{\lra}[s_{i-1}]
\stackrel{g_{i-1}}{\lra}\ldots \stackrel{g_{2}}{\lra}[s_1]} F(t),
\end{equation}

whose homology is denoted by
$H_{(s_n,s_{n-1},\ldots,s_{i+1},*,s_{i-1},\ldots,s_1)}(F,\partial_i)$.
There is an analogous splitting for the complex
$(C_*^{E_n,X}(F),\partial_i)$.

In the proposition below we compute $H_*(F,\partial_2)$ for the
representable  functors for $n=2$, give its generators
in propositions \ref{prop:tophomology} and \ref{prop:shufflecycles}
and then  prove that the representable functors are acyclic for $n=2$ in
\ref{cor:nonfork2} and \ref{cor:fork2} depending on the form of the
$2$-level  tree $t$, whether it is a fork tree or not.

\begin{prop}\label{prop:Hd2} Let $(t,\phi)=
  X\stackrel{\phi}{\rightarrow}[r_2]\stackrel{f}{\lra} [r_1]$ be
  an $(X,2)$-level tree in $\epi_2^X$.

$$\begin{array}{ll}
H_{(*,s)}({\epi_2^X}^{t,\phi},\partial_2)=\ \ 0 & \hbox{\rm if\ } r_2\not=r_1, \\
H_{(*,s)}({\epi_2^X}^{t,\phi},\partial_2)\cong \begin{cases} 0& \hbox{\rm for}\
  *\not=r_2, \\
k^{\oplus |\deltaepi([r_2],[s])|}& \hbox{\rm for}\  s\leq *=r_2,
\end{cases} &    \hbox{\rm if\ } r_2=r_1.\\
\end{array}$$
\end{prop}

\begin{proof}
Let $F$ denote the covariant functor ${\epi_2^X}^{t,\phi}$.

Assume $s=0$. We first prove that the chain complex $\partial_2\colon
C_{(*,0)}^{E_2,X}(F)\rightarrow C_{(*-1,0)}^{E_2,X}(F)$ is the normalized chain
complex associated to a small category as in definition \ref{def:dgcat}.

The chain complex $(C_{(*,0)}^{E_2,X}(F),\partial_2)$ has the following 
form, for $0 < u \leq r_2$

$$ \xymatrix@1{ {\bigoplus_{\psi} k[\epi_2^X((t,\phi);
X\stackrel{\psi}{\rightarrow}[u]\stackrel{}{\rightarrow} [0])]}
\ar[rrr]^(0.45){\sum_{i=0}^{u-1}(-1)^{s_{2,i}} d_i^2} & & &
{\bigoplus_{\psi} k[\epi_2^X((t,\phi);X
\stackrel{\psi}{\rightarrow}[u-1]
\stackrel{}{\rightarrow}[0])].}}$$

Let $(A_0,\ldots,A_{r_1})$ be the sequence of preimages of $f$, and
$a_i$ the number of elements in $A_i$. For a fixed $\psi$, the set
$\epi_2^X((t,\phi); 
X\stackrel{\psi}\rightarrow[u]{\rightarrow}[0])$ is either empty or
has only one element uniquely determined by a surjective map
$\sigma\colon[r_2]\rightarrow [u]$: 
since $\phi$ is surjective,
the requirement  $\psi=\sigma\phi$ uniquely determines the surjection 
$\sigma$ if it exists. In that case, $\sigma$ determines an element in 
$\epi_2^X((t,\phi); 
X\stackrel{\psi}\rightarrow[u]{\rightarrow}[0])$ if it is
order-preserving on the fibres of $f$. 

The map
$\sigma$ can be described by the sequence of its preimages
$(S_0,\ldots, S_u)$ with the condition $(C_S)$: if $a<b\in A_i$ then
$i_a\leq i_b$ where $i_\alpha$ is the unique index for which
$\alpha\in S_{i_\alpha}$.
One has
$$d(S_0,\ldots,S_u)=\sum_{i=0}^{u-1} (-1)^{i+d(\phi^{-1}(S_0\cup\ldots\cup 
S_i))}\varepsilon(\phi^{-1}(S_i);\phi^{-1}(S_{i+1}))(S_0,\ldots,S_i\cup 
S_{i+1},\ldots,S_u).$$ Let $\mathcal C$ be the tensor product of the
categories $[a_j]^{\pi_j}, 0\leq j\leq r_1$ where the 
partitions $\pi_j$'s will be
defined later. The category $\mathcal C$ is a poset category and the
order is given by the product order. Any object of the category
$\mathcal C$ can be written as $p=(p^0,\ldots,p^{r_1})$. We denote
by $0$ the minimal element $(0,\ldots,0)$, by $\alpha$ the maximal 
element $(a_0,\ldots,a_{r_1})$ and 
by $L_0$ and $R_\alpha$ the corresponding
left and right graded $\mathcal C$-module. For $p\leq q$ the element
$\xi_{pq}$ denotes the unit of $k=\mathcal C(p,q)$. The reduced
complex (\ref{eq:normalized}) computing $\mathrm{Tor}^{\mathcal
C}_{u+1}(R_\alpha;L_0)$ is given by
$$C_{u+1}(\mathcal C)=L_0\otimes_{\mathcal{C}^\mathrm{op}} N_{u+1}(\mathcal{C})
\otimes_{\mathcal{C}^\mathrm{op}}
R_\alpha=\bigoplus_{0<p_1<\ldots<p_u<\alpha} k\
\xi_{0p_1}\otimes\xi_{p_1p_2}\ldots\otimes\xi_{p_u\alpha},$$ with
the differential $d=\sum_{i=1}^{u}(-1)^i d_i$ given in definition
\ref{def:dgcat}. The sequence $(S_0,\ldots,S_u)$ is in 1-to-1
correspondence with the element $
\xi_{0p_1}\otimes\xi_{p_1p_2}\ldots\otimes\xi_{p_u\alpha}$ where the
$j$-th coordinate of $p_i$, $p_i^j$, is given by the number of elements in
$A_j\cap (\cup_{k\leq i-1} S_k)$. Conversely, let $p_0=0<p_1<\ldots 
<p_u<p_{u+1}=\alpha$ be a sequence of objects in $\otimes [a_j]^{\pi_j}$. 
The set $A_i$ is  
$\{\sum_{k=0}^{i-1} a_k +l, 0\leq l\leq a_{i}-1\}$. Let $S_i\subset [r_2]$ 
be the set defined by $S_i\cap A_j=\{\sum_{k=0}^{j-1} a_k +l,  
p_i^j\leq l\leq p_{i+1}^j-1\}$. The sequence $(S_0,\ldots,S_u)$ satisfies 
the condition $(C_S)$. The partition $\pi_j$ is
the partition of $(f\phi)^{-1}(j)$ given by 
$\pi_j^k=\phi^{-1}(k)$ whenever $f(k)=j$. As a consequence, the two
complexes, $(C_{u+1}(\mathcal C),-d)$ and
$(C_{(u,0)}^{E_2,X}(F),\partial_2)$, coincide. By lemma
\ref{lem:torcomputation} the homology
$H_{(u,s)}({\epi_2^X}^{t,\phi},\partial_2)$ is $0$ but for the case where 
$u+1=r_1+1$ and $a_j=1$ for all $j$. The latter condition is
equivalent to $r_2=r_1$ and $f=\id$. This concludes the case $s=0$.

\medskip

Assume $s>0$.

The complex $(C_{(*,s)}^{E_2,X}(F),\partial_2)$ splits into subcomplexes
$$C_{(*,s)}^{E_2,X}(F)=\bigoplus_{\sigma\in\deltaepi([r_1],[s])} C_{(*,s)}(F_\sigma)=
\bigoplus_{\sigma\in\deltaepi([r_1],[s])}\bigoplus_{g\in\deltaepi([*],[s]),\psi}
F_\sigma(X\stackrel{\psi}{\lra}[*] 
\stackrel{g}{\lra}[s])$$
where $F_\sigma(X\stackrel{\psi}{\lra}[u]\stackrel{g}{\lra}[s])\subset
{\epi_2^X}^{t,\phi}(X\stackrel{\psi}{\lra}[u]
\stackrel{g}{\lra}[s])$
is the free $k$-module generated by morphisms of the form
\begin{equation}\label{eq:notation}
 \xymatrix{X\ar[r]^{\phi}\ar[d]^{\id}&[r_2] \ar[d]^{\tau} \ar[r]^{f} &{[r_1]}
  \ar[d]^{\sigma}\\
 X\ar[r]^{\psi}&{[u]} \ar[r]^{g} & {[s].}}
\end{equation}

Let $(A_0,\ldots,A_s)$ denote the sequence of preimages of $\sigma
f$ and $(B_0,\ldots,B_s)$ that of $g$. The latter has to satisfy the
condition $|B_i|\leq |A_i|,0\leq i\leq s$.  Note that
$g\in\deltaepi([u],[s])$ is also  uniquely determined by the
sequence $(b_0,\ldots,b_s)$ of the cardinalities of its preimages.
The differential $\partial_2\colon C_{(u,s)}(F_\sigma)\lra
C_{(u-1,s)}(F_\sigma)$ has the following form:

\begin{multline*}
\partial_2\left(\vcenter{\xymatrix{X\ar[r]^{\phi}\ar[d]^{\id}&[r_2]\ar[d]^{\tau}
\ar[r]^{f}  & [r_1]
  \ar[d]^{\sigma}\\
 X\ar[r]^{\tau\phi}&[u] \ar[r]^{g} & [s]}}\right)
 = \sum_{i|g(i)=g(i+1)}(-1)^{s_{2,i}}\epsilon((\tau\phi)^{-1}(i);
 (\tau\phi)^{-1}(i+1))   
\vcenter{\xymatrix{X\ar[r]^{\phi}\ar[d]^{\id}&[r_2]\ar[d]^{d_i\tau}
    \ar[r]^{f} &  [r_1]
  \ar[d]^{\sigma}\\
 X\ar[r]^{d_i\tau\phi}&[u-1] \ar[r]^{g|_{i=i+1}} & [s]}} \\
 =\sum_{j=0}^s\left(\sum_{i\in B_j|g(i)=g(i+1)} 
(-1)^{s_{2,i}}\epsilon((\tau\phi)^{-1}(i); (\tau\phi)^{-1}(i+1)) 
\vcenter{\xymatrix{X\ar[r]^{\phi}\ar[d]^{\id}&[r_2]\ar[d]^{d_i\tau}
    \ar[r]^{f} &  [r_1]
  \ar[d]^{\sigma}\\
 X\ar[r]^{d_i\tau\phi}&[u-1] \ar[r]^{g|_{i=i+1}} & [s]}}\right).
\end{multline*}

Define $D_j$ by
restricting the sum over indices $i$ such that $g(i)=g(i+1)$ to the
sum over  indices $i\in B_j$ such that $g(i)=g(i+1)$. One has
$$D_j\colon C_{(u,s)}(F_\sigma)=\bigoplus_{b_0+\ldots+b_s=u+1}
C_{((b_0,\ldots,b_j,\ldots,b_s),s)}(F_\sigma) \lra
\bigoplus_{b_0+\ldots+b_s=u+1}
C_{((b_0,\ldots,b_j-1,\ldots,b_s),s)}(F_\sigma)$$ and
$\partial_2=D_0+\ldots+D_s$. We claim that the $D_j$ are
anti-commuting differentials:

Let $i$ be in $B_j$ and $\ell$ be in $B_k$. For $j<k$ it follows that
$i+1<\ell$ for $g(i)=g(i+1)=j<g(l)=k$ therefore we have the relation
$d_i  d_\ell = d_{\ell-1} d_i$. Furthermore the $\epsilon$-signs
involved do not depend on the way we  compose: the resulting sign is
$\epsilon((\tau\phi)^{-1}(i);(\tau\phi)^{-1}(i+1))\epsilon((\tau\phi)^{-1}(l);
(\tau\phi)^{-1}(l+1))$ for $d_{i}^{-1}(l-1)=l$ and $d_l^{-1}(i)=i$.
In order to calculate the effect of $d_{\ell-1}
d_i$ we have to determine $s_{2,\ell-1}$ after the application of
$d_i$. Let $\tilde{S}_j$ denote the preimage $(d_i \circ
\tau\circ\phi)^{-1}(j)$ and $S_j$ the preimage
$(\tau\circ\phi)^{-1}(j)$ for $j \in [u]$.
Then
$$ d(\tilde{S}_j) = \begin{cases}  d(S_j), & j < i\\
d(S_i)+d(S_{i+1}), & j=i \\
d(S_{j+1}), & j > i. \end{cases}$$
Thus $s_{2,\ell-1}$ is $\ell
-1+k+2 +\sum_{j=0}^{\ell-1} d(\tilde{S}_j) = \ell+k +1 +
\sum_{j=0}^{\ell} d(S_j)$ whereas $s_{2,\ell} = \ell+k +2 +
\sum_{j=0}^{\ell} d(S_j)$. A similar argument shows that the $D_j$
are differentials.

The complex $(C_{(u,s)}(F_\sigma), D_s)$ splits  into subcomplexes
$(C_{((b_0,\ldots,b_{s-1}),*)}(F_\sigma),D_s)$ for fixed $b_i\leq
a_i=|A_i|, i<s$. With the notation of definition
\ref{def:enchainsX}, the tree $(t,\phi)$  can be written as
$t=[t_{1,0},\ldots,t_{1,r_1}]$, with $t_{1,i}$ being an
$(X_{1,i},1)$-level tree. Let $p$ be the first integer such that
$\sigma(p)=s$. Let $X_{s-1}=\cup_{0\leq
  i\leq p-1}X_{1,i}$ and
$\tilde X=\cup_{p\leq i\leq r_1}X_{1,i}$. Denote by $t_{s-1}$ the
$(X_{s-1},2)$-level tree $t_{s-1}=[t_{1,0},\ldots,t_{1,p-1}]$ and by
$\tilde t$ the $(\tilde X,2)$-level tree $\tilde
t=[t_{1,p},\ldots,t_{1,r_1}]$. See \ref{ex:cuttree} for an example.
Let $\sigma_{s-1}$ (resp.
$\phi_{s-1}$) be the map obtained from $\sigma$ (resp. $\phi$) by
restriction
$\sigma_{s-1}\colon\sigma^{-1}([s-1])\stackrel{\sigma}{\lra} [s-1]$.
Let $u_{s-1}=(\sum_{i<s}b_i)-1$. The subcomplex
$(C_{((b_0,\ldots,b_{s-1}),*)}(F_\sigma),D_s)$ can be expressed as
$$\bigoplus_{\psi;\gamma\in({\epi_2^{X_{s-1}}}^{t_{s-1},\phi_{s-1}})_{\sigma_{s-1}}
  (X_{s-1}\stackrel{\psi}{\lra}[u_{s-1}] 
\stackrel{g|_{[u_{s-1}]}}{\lra}[s-1])} (C^{E_2,\tilde
X}_{(*,0)}({\epi_2^{\tilde X}}^{\tilde t,\tilde\phi}),
(-1)^{b_0+\ldots+b_{s-1}+s+d(X_{s-1})}\partial_2).$$

If $f\not=\id$, then  there exists $j\in[s]$ such that the
restriction of $f$ on $(\sigma\circ f)^{-1}(j)\ra \sigma^{-1}(j)$ is
different from the identity.  Without loss of generality we can
assume that $j=s$, hence $\tilde  t$ is a non-fork tree and the
homology of the complex is $0$. If  $f=\id$, then we deduce from the
case $s=0$  that the complex $(C^{E_2,\tilde
X}_{(*,0)}({\epi_2^{\tilde X}}^{\tilde t,\tilde\phi}),\partial_2)$
has only top homology of  rank one; consequently when
$t\colon[r_2]\lra[r_2]$ is the fork tree,
$$(H_*(C_{(*,s)}(({\epi_2^X}^{t,\phi})_\sigma),D_s),D_1+\ldots+D_{s-1})\cong
(C_{(*,s-1)}(({\epi_2^{X_{s-1}}}^{t_{s-1},\phi_{s-1}})_{\sigma_{s-1}}),\partial_2).$$
We then have an inductive process to compute the homology of the total complex
$(C_{(*,s)}(F_\sigma),\partial_2)$. Consequently, for a fixed
$\sigma\colon[r_2] \ra[s]$
$$\begin{array}{ll}
H_{(*,s)}(F_\sigma,\partial_2)= 0, &   \hbox{\rm if\ } r_2\not=r_1 \\
H_{(*,s)}(F_\sigma,\partial_2)\cong \begin{cases}  0& \hbox{\rm  for}\
  *\not= r_2 \\
k& \hbox{\rm for}\  s\leq *=r_2 \end{cases} ,&   \hbox{\rm if\ } r_2=r_1. \\
\end{array}$$
Since each $\sigma\in\deltaepi([r_2],[s])$ contributes to one summand in
$H_{(r_2,s)}(F,\partial_2)$, this proves the claim. .
The computation of the generators for $s>0$ is given in proposition
\ref{prop:shufflecycles}.
\end{proof}

\begin{ex} \label{ex:cuttree}
Let $(t,\phi)=X\stackrel{\phi}{\lra}[6]\stackrel{f}{\lra}[2]$ be the
following tree

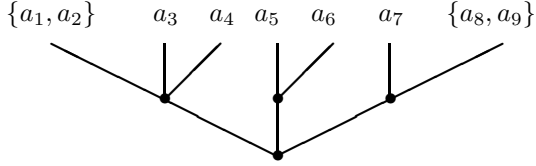
\raisebox{-3cm}{
\begin{picture}(10,10)
\setlength{\unitlength}{1.5cm} \thicklines
\put(3,0.5){\line(-2,1){2}} \put(3,0.5){\line(0,1){1}}
\put(3,0.5){\line(2,1){2}} \put(3,1){\line(1,1){0.5}}
\put(2,1){\line(1,1){0.5}} \put(2,1){\line(0,1){0.5}}
\put(4,1){\line(0,1){0.5}} \put(0.6,1.7){$\{a_1,a_2\}$}
\put(1.9,1.7){$a_3$} \put(2.4,1.7){$a_4$} \put(2.8,1.7){$a_5$}
\put(3.3,1.7){$a_6$} \put(3.9,1.7){$a_7$}
\put(4.5,1.7){$\{a_8,a_9\}$} \put(2.95,0.45){$\bullet$}
\put(1.95,0.95){$\bullet$} \put(2.95,0.95){$\bullet$}
\put(3.95,0.95){$\bullet$}
\end{picture}}

where $(t,\phi)=[t_{1,0},t_{1,1},t_{1,2}]$ with

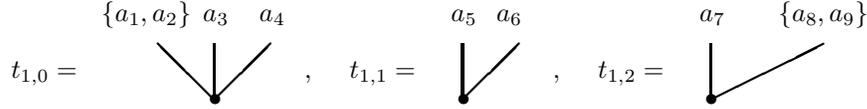
\raisebox{-3cm}{
\begin{picture}(10,10)
\setlength{\unitlength}{1.5cm} \thicklines \put(0.2,1.2){$t_{1,0}=$}
\put(2,1){\line(1,1){0.5}} \put(2,1){\line(0,1){0.5}}
\put(2,1){\line(-1,1){0.5}} \put(1.0,1.7){$\{a_1,a_2\}$}
\put(1.9,1.7){$a_3$} \put(2.4,1.7){$a_4$} \put(1.95,0.95){$\bullet$}
\put(2.8,1.2){,} \put(3.2,1.2){$t_{1,1}=$}
\put(4.2,1){\line(0,1){0,5}} \put(4.2,1){\line(1,1){0,5}}
\put(4.15,0.95){$\bullet$} \put(4.1,1.7){$a_5$} \put(4.5,1.7){$a_6$}
\put(5,1.2){,} \put(5.4,1.2){$t_{1,2}=$}
\put(6.4,1){\line(0,1){0,5}} \put(6.4,1){\line(2,1){1}}
\put(6.35,0.95){$\bullet$} \put(6.3,1.7){$a_7$}
\put(7,1.7){$\{a_8,a_9\}$}
\end{picture}}

and $X_{1,0}=\{a_1,a_2,a_3,a_4\}, X_{1,1}=\{a_5,a_6\}$ and
$X_{1,2}=\{a_7,a_8,a_9\}$. Let $\sigma:[2]\rightarrow [1]$ be the
map assigning $0$ to $0$ and $1$ to $1$ and $2$. One has $s=1$,
$p=1$, so that $X_{s-1}=\{a_1,\ldots,a_4\}$ and $\tilde
X=\{a_5,\ldots, a_9\}$. Moreover

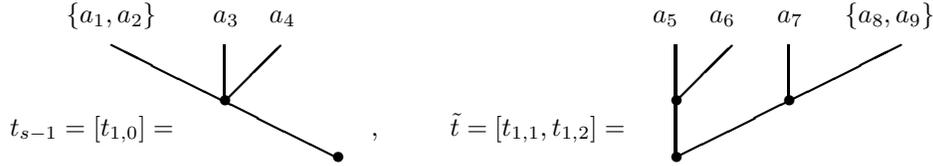
\raisebox{-3cm}{
\begin{picture}(10,10)
\setlength{\unitlength}{1.5cm} \thicklines
\put(0.1,0.7){$t_{s-1}=[t_{1,0}]=$} \put(3,0.5){\line(-2,1){2}}
\put(2,1){\line(1,1){0.5}} \put(2,1){\line(0,1){0.5}}
\put(0.6,1.7){$\{a_1,a_2\}$} \put(1.9,1.7){$a_3$}
\put(2.4,1.7){$a_4$}
\put(2.95,0.45){$\bullet$} \put(1.95,0.95){$\bullet$}
\put(3.3,0.7){,} \put(4,0.7){$\tilde t=[t_{1,1},t_{1,2}]=$}
\put(6,0.5){\line(0,1){1}} \put(6,0.5){\line(2,1){2}}
\put(6,1){\line(1,1){0.5}} \put(7,1){\line(0,1){0.5}}
\put(5.8,1.7){$a_5$} \put(6.3,1.7){$a_6$} \put(6.9,1.7){$a_7$}
\put(7.5,1.7){$\{a_8,a_9\}$} \put(5.95,0.45){$\bullet$}
\put(5.95,0.95){$\bullet$} \put(6.95,0.95){$\bullet$}
\end{picture}}

\end{ex}

\begin{cor}\label{cor:nonfork2} For any non-fork tree $(t,\phi)=
  X\stackrel{\phi}{\lra}[r_2] 
\stackrel{f}{\lra} [r_1], r_2\not=r_1$, ${\epi_2^X}^{t,\phi}$ is acyclic.
For any non-fork tree $t=[r_2]
\stackrel{f}{\lra} [r_1], r_2\not=r_1$, $\epi_2^t$ is acyclic.
\end{cor}
\begin{proof}
The first assertion is a direct consequence of the first equation of
proposition  \ref{prop:Hd2}. The second one is a direct consequence of 
relation (\ref{eq:Xnongraded}).
\end{proof}

\begin{prop}\label{prop:tophomology}
Let $(t,\phi)\colon
X\stackrel{\phi}{\lra}[r]\stackrel{\id}{\lra}[r]$ be a fork tree and
let $X_i=\phi^{-1}(i)$. Then the top homology
$H_{(r,0)}({\epi_2^X}^{t,\phi},\partial_2)$ is freely generated by
$c_{r,X} := \sum_{\sigma \in \Sigma_{r+1}} \mathrm{sgn}(\sigma;X)
\sigma$, where the sign $\mathrm{sgn}(\sigma;X)$ picks up a factor
$(-1)^{(d(X_i)+1)(d(X_j)+1)}$ whenever
$\sigma(i)>\sigma(j)$ but $i<j$.

In particular, for a fork tree $t\colon[r]\stackrel{\id}{\lra}[r]$,
the top homology
$H_{(r,0)}({\epi_2}^t,\partial_2)$ is freely generated by
$c_{r} := \sum_{\sigma \in \Sigma_{r+1}} \mathrm{sgn}(\sigma) \sigma$.

\end{prop}

\begin{proof}
The second assertion is a consequence of the first one using relation
(\ref{eq:Xnongraded}).
The computation of the top homology amounts to determining
  the  kernel of the map
$$\xymatrix@1{
\bigoplus_{\psi}k[\epi_2^X(X\stackrel{\phi}{\lra}[r]\stackrel{\id}{\lra}[r];X
\stackrel{\psi}{\lra}[r]\lra [0])] \ar[r]^{\partial_2} &
\bigoplus_{\psi}
k[\epi_2^X(X\stackrel{\phi}{\lra}[r]\stackrel{\id}{\lra}[r]; X
\stackrel{\psi}{\lra}[r-1]\lra [0])].}$$
The set
$\epi_2^X(X\stackrel{\phi}{\lra}[r]\stackrel{\id}{\lra}[r];X
\stackrel{\psi}{\lra}[r]\lra [0])$ is either empty or has only one
element uniquely determined by the following diagram
$$\xymatrix{X\ar[r]^{\phi}\ar[d]_{\id}&[r]\ar[d]^{\tau}\ar[r]^{\id} & [r]\ar[d]
  \\
 X\ar[r]^{\psi=\tau\phi}&[r]\ar[r]  & [0].}$$
Hence, the surjection $\psi$ determines a bijection
$\tau$  and this induces a permutation of the set
$\{X_0,\ldots,X_r\}$.  We denote such an element by $\tau \cdot X := 
(X_{\tau^{-1}(0)},\ldots,X_{\tau^{-1}(r)})$. As a consequence the
computation of the top homology amounts to determining
  the  kernel of the map
$$\partial_2\colon k[\Sigma_{r+1}] \lra \bigoplus_{\psi}
k[\epi_2^X(X\stackrel{\phi}{\lra}[r]\stackrel{\id}{\lra}[r];X
\stackrel{\psi}{\lra}[r-1]\lra [0])] $$
where
$$\partial_2(\tau\cdot X)=\sum_{i=0}^{r-1}
(-1)^{i+d(X_{\tau^{-1}(0)})+\ldots+d(X_{\tau^{-1}(i)})}\epsilon(X_{\tau^{-1}(i)};
X_{\tau^{-1}(i+1)})
(X_{\tau^{-1}(0)},\ldots, X_{\tau^{-1}(i)}\cup
X_{\tau^{-1}(i+1)},\ldots,X_{\tau^{-1}(r)} ).$$ Therefore, if
$x=\sum_{\tau\in\Sigma_{r+1} }\lambda_\tau \tau \cdot X$ is in the
kernel of $\partial_2$, then for all transpositions $(i,i+1)$ and
all $\tau$ one has 
$\lambda_{(i,i+1)\tau}=(-1)^{1+d(X_{\tau^{-1}(i)})+d(X_{\tau^{-1}(i+1)})+
d(X_{\tau^{-1}(i)})d(X_{\tau^{-1}(i+1)})}\lambda_\tau$. Since the
transpositions generate the symmetric group one has $\lambda_\tau=
\mathrm{sgn}(\tau;X) \lambda_{\id}$ and $x=\lambda_{\id} c_{r,X}$.
\end{proof}

For $s>0$, the computation of the top homology of
$(C_{*,s}^{E_2,X}({\epi_2^X}^{t,\phi}),\partial_2)$ amounts to
calculating
the  kernel of the map $\partial_2$
$$\bigoplus_{\psi,g\in\deltaepi([r],[s])}
k[\epi_2^X((t,\phi);X\stackrel{\psi}{\lra}[r]\stackrel{g}{\lra} [s])] \lra
\bigoplus_{\psi,h\in\deltaepi([r-1],[s])}k[\epi_2^X((t,\phi);X\stackrel{\psi}{\lra}
[r-1] \stackrel{h}{\lra} [s])].$$
We know from proposition \ref{prop:Hd2} that it is free of rank equal to the
cardinality of  $\deltaepi([r],[s])$.
As before, the set
$\epi_2^X((t,\phi);X\stackrel{\psi}{\lra}[r]\stackrel{g}{\lra} [s])$
is either empty or has only one element determined by the commuting diagram

$$\xymatrix{X\ar[r]^{\phi}\ar[d]^{\id}&[r]\ar[d]^{\tau}\ar[r]^{\id} &
  [r]\ar[d]^{g'}\\ 
 X\ar[r]^{\psi=\tau\phi}&[r]\ar[r]^{g}&[s]  .}$$

An element $g$ in $\deltaepi([r],[s])$ is uniquely determined by the
sequence  $(x_0,\ldots,x_s)$ of the cardinalities of its preimages. 
Furthermore,  any map in $\epi_2([r]\stackrel{\id}{\lra}[r];[r]
\stackrel{g}{\lra} [s])] $ is given by $g' \colon[r]\ra[s]$ in
$\deltaepi$ and $\tau\colon[r]\ra [r]$ in $\Sigma_{r+1}$  such that
$g' =g\tau$. This implies that $g'=g$ and $\tau\in
\Sigma_{x_0}\times\ldots\times \Sigma_{x_s}$. If there is such a
$\tau$ satisfying $\psi=\tau\phi$ then the set is non-empty and
$\tau$ is unique. Let $X_{(x_i)}=(g\phi)^{-1}(i)$. Then $X_{(x_i)}$
is a subset of $X$ and there is a natural partition of it given by
$X_{(x_i)}=\sqcup_{j \in g^{-1}(\{i\})} X_j$.

Let $c_{(x_0,\ldots,x_s);X}$ be the element
$$c_{(x_0, \ldots,
x_s);X} = (\sum_{\sigma^0 \in
\Sigma_{x_0}}\mathrm{sgn}(\sigma^0;X_{(x_0)})\sigma^0,\ldots
,\sum_{\sigma^s \in
\Sigma_{x_s}}\mathrm{sgn}(\sigma^s;X_{(x_s)})\sigma^s) \in
k[\Sigma_{x_0}\times\ldots\times\Sigma_{x_s}].$$ If every element of
$X$ has degree zero, we denote  $(\sum_{\sigma^0 \in
\Sigma_{x_0}}\mathrm{sgn}(\sigma^0)\sigma^0,\ldots ,\sum_{\sigma^s
\in \Sigma_{x_s}}\mathrm{sgn}(\sigma^s)\sigma^s)$ by
$c_{(x_0,\ldots,x_s)}$.

\begin{prop} \label{prop:shufflecycles}
Let $(t,\phi)\colon X\stackrel{\phi}{\lra}[r]\stackrel{\id}{\lra}[r]$
be a fork tree. The top homology $H_{(r,s)}({\epi_2^X}^{t,\phi},\partial_2)$ is
freely generated by the elements $c_{(x_0, \ldots,x_s);X} =
(\sum_{\sigma^0 \in
\Sigma_{x_0}}\mathrm{sgn}(\sigma^0;X_{(x_0)})\sigma^0,\ldots
,\sum_{\sigma^s \in
\Sigma_{x_s}}\mathrm{sgn}(\sigma^s;X_{(x_s)})\sigma^s)$, for
$g=(x_0,\ldots,x_s)\in  \deltaepi([r],[s])$,
$X_{(x_k)}=(g\phi)^{-1}(\{k\})$.

Let $t\colon[r]\stackrel{\id}{\lra}[r]$ be a fork tree. The top
homology $H_{(r,s)}(\epi_2^t,\partial_2)$ is freely generated by the
elements $c_{(x_0, \ldots,x_s)} = (\sum_{\sigma^0 \in
\Sigma_{x_0}}\mathrm{sgn}(\sigma^0)\sigma^0,\ldots ,\sum_{\sigma^s
\in \Sigma_{x_s}}\mathrm{sgn}(\sigma^s)\sigma^s)$, for
$(x_0,\ldots,x_s)\in \deltaepi([r],[s])$.
\end{prop}
\begin{proof}
As in the proof of proposition \ref{prop:tophomology} we compute the
kernel of $\partial_2$ which decomposes into the sum of
anti-commuting differentials $\partial_2=D_0+\ldots+D_s$, as in the
proof of proposition \ref{prop:Hd2}. As a consequence
$\ker(\partial_2)=\cap_i \ker(D_i)$, which  gives the result.
\end{proof}

\begin{cor}\label{cor:fork2} For any fork tree
  $(t,\phi)=X\stackrel{\phi}{\lra}[r]\stackrel{\id}{\lra}  [r]$,
${\epi_2^X}^{t,\phi}$ is acyclic. In particular, $\epi_2^{t}$ is
acyclic for any fork tree $t=  [r]\stackrel{\id}{\lra}  [r]$.
\end{cor}
\begin{proof}
It remains to compute the homology of the complex
$((H_{(r,*)}(C^{E_2,X}({\epi_2^X}^{t,\phi}),\partial_2),\partial_1)$
and prove that it vanishes for all $*$ if $r>0$. From propositions
\ref{prop:tophomology} and
\ref{prop:shufflecycles} one has
$$H_{(r,s)}(C^{E_2,X}({\epi_2^X}^{t,\phi}),\partial_2)=\bigoplus_{(x_0,\ldots,x_s)
  \in    \deltaepi([r],[s])} kc_{(x_0,\ldots,x_s);X}.$$
To compute $\partial_1(c_{(x_0,\ldots,x_s);X})$ it is enough to
compute $\partial_1(\id_{\Sigma_0\times \ldots\times\Sigma_s})$ in
$C^{E_2,X}_{(r,s)}({\epi_2^X}^{t,\phi})$. We apply relations
(\ref{eq:dijX}) and (\ref{eq:cenX}):

\begin{multline*}
\partial_1
\left(\vcenter{\xymatrix{X\ar[r]^{\phi}\ar[d]^{\id}&[r]\ar[r]^{\id}\ar[d]^{\id}&
      [r] 
\ar[d]^{(x_0,\ldots,x_s)}\\
X\ar[r]^{\phi}&[r]\ar[r]^{(x_0,\ldots,x_s)}&[s]}}\right)
 \\
=
\sum_{i=0}^{s-1}(-1)^{i+1+x_0+d(X_{(x_0)})+\ldots+x_i+d(X_{(x_i)})}\left(\vcenter{
    \xymatrix{ 
X\ar[r]^{\phi}\ar[d]^{\id}&[r]\ar[r]^{\id}\ar[d]^{\id}&
      [r]\ar[d]^{d_i(x_0,\ldots,x_s)}\\
   X\ar[r]^{\phi}&   [r]\ar[r]^{d_i(x_0,\ldots,x_s)}&[s-1]}}\right.
   \\
   \pm  \left.\sum_{\xi}
\vcenter{\xymatrix{X\ar[r]^{\phi}\ar[d]^{\id}&[r]\ar[r]^{\id}\ar[d]^{\xi}&
    [r]\ar[d]^{d_i(x_0,\ldots,x_s)} \\
 X\ar[r]^{\xi\phi}&   [r]\ar[r]^{d_i(x_0,\ldots,x_s)}&[s-1]}} \right),
\end{multline*}
with $\xi$ running over the $(X_{(x_i)},X_{(x_{i+1})})$-shuffles with
$\xi\not=\id$.
Thus,
$$\partial_1(c_{(x_0,\ldots,x_s);X}
)=\sum_{i=0}^{s-1}(-1)^{i+1+x_0+d(X_{(x_0)})+ \ldots+x_i+d(X_{(x_i)})}
c_{(x_0,\ldots,x_i+x_{i+1}, \ldots,x_s);X}$$ and the complex
$(H_{(r,*)}(C^{E_2,X}({\epi_2^X}^{t,\phi}),\partial_2),\partial_1)$
agrees with the graded version of the  complex
$C_*^{\mathrm{bar}}(({\deltaepi})^r)$ of remark \ref{R:gradedn=1}.
Therefore it is acyclic, with
$$H_0(C_*^{\mathrm{bar}}(({\deltaepi})^r)=\begin{cases} 0 & \mathrm{if}\
  r>0 \\ k&  \mathrm{if}\ r=0.\end{cases}$$
As a consequence the spectral sequence associated to the bicomplex
$(C^{E_2,X}_{(*,*)}({\epi_2^X}^{t,\phi}),\partial_1+\partial_2)$ collapses at the 
$E^2$-stage and one gets $H_p^{E_2,X}({\epi_2^X}^{t,\phi})=0$ for all $p>0$.
\end{proof}

\begin{prop}\label{prop:Hdn} Let $X=\{x_0<\ldots<x_{r_n}\}$ be an ordered set 
of graded elements. Let
$\phi:X\rightarrow [r_n]$ be the map sending $x_i$ to $i$. Let
$(t,\phi)=X\stackrel{\phi}{\lra}[r_n]\stackrel{f_n}{\lra}[r_{n-1}] 
\stackrel{f_{n-1}}{\lra}\ldots 
\stackrel{f_2}{\lra}[r_1]$ be an $(X,n)$-level tree and let $\bar t$
be its $(n-1)$-truncation
$X[1]\stackrel{f_{n}\phi}{\lra}[r_{n-1}]\stackrel{f_{n-1}}{\lra}\ldots
\stackrel{f_2}{\lra}[r_1]$,
where $X[1]$ is the ordered set obtained from $X$ by increasing the
degree of its elements by $1$, then
$$\begin{array}{ll}
H_{(*,s_{n-1},\ldots,s_1)}({\epi_n^X}^{t,\phi},\partial_n)= 0, &
\hbox{\rm if\
} r_{n} \not=r_{n-1}, \\
H_{(*,s_{n-1},\ldots,s_1)}({\epi_n^X}^{t,\phi},\partial_n)\cong
\begin{cases}  0&
  \hbox{\rm  for}\   *\not=r_n \\
C_{(s_{n-1},\ldots,s_1)}^{E_{n-1},X[1]}({\epi_{n-1}^{X[1]}}^{{\bar
    t},\phi})&  \hbox{\rm
  for}\  s_{n-1}\leq  *=r_n \end{cases} ,&   \hbox{\rm if\ } r_n=r_{n-1}. \\
\end{array}$$
Furthermore when $f_n=\id$, the $(n-1)$-complex structure induced on
$H_{(r_n,s_{n-1},\ldots,s_1)}({\epi_n^X}^{t,\phi},\partial_n)$ by
the $n$-complex structure of
$C_{(*,\ldots,*)}^{E_n,X}({\epi_n^X}^{t,\phi})$ coincides with the
one on
$C_{(s_{n-1},\ldots,s_1)}^{E_{n-1},X[1]}({\epi_{n-1}^{X[1]}}^{\bar
t,\phi})$.
\end{prop}
\begin{proof}
Recall from definition \ref{def:enchainsX} that
\begin{multline*}
\partial_n\left(\vcenter{\xymatrix{X\ar[r]^{\phi}\ar[d]^{\id}&[r_n]\ar[r]^{f_n}
      \ar[d]^{\sigma_n}&
[r_{n-1}]\ar[r]^{f_{n-1}}\ar[d]^{\sigma_{n-1}}&
\ldots\ar[r]^{f_2}&[r_1]\ar[d]^{\sigma_1}\\
X\ar[r]^{\sigma_n\phi}&[s_n]\ar[r]^{g_n}& [s_{n-1}]\ar[r]^{g_{n-1}}&
\ldots\ar[r]^{g_2}&[s_1]}}\right) \\
=
\sum_{i|g_n(i)=g_n(i+1)}(-1)^{s_{n,i}}\epsilon((\sigma_n\phi)^{-1}(i);
(\sigma_n\phi)^{-1}(i+1)) 
\vcenter{\xymatrix{X\ar[r]^{\phi}\ar[d]^{\id}&[r_n]\ar[r]^{f_n}\ar[d]^{d_i\sigma_n}
 & [r_{n-1}]\ar[r]^{f_{n-1}} \ar[d]^{\sigma_{n-1}}&
\ldots\ar[r]^{f_2}&[r_1]\ar[d]^{\sigma_1}\\
X\ar[r]^(0.4){d_i\sigma_n\phi}&[s_n-1]\ar[r]^{g_n|_{i=i+1}}&
[s_{n-1}]\ar[r]^{g_{n-1}}& \ldots\ar[r]^{g_2}&[s_1].}}
\end{multline*}
The same proof as in proposition \ref{prop:Hd2} provides the
computation of  the homology of the complex
with respect to the differential
$\partial_n$: if $t$ is not a fork tree, then the homology of the
complex vanishes, and if $t$ is  the fork tree $f_n=\id_{[r_{n-1}]}$,
then its homology groups are concentrated in top degree $r_n$.  Let us
describe all the bijections $\tau$ of
$[r_{n-1}]$ such that the following diagram commutes
$$\xymatrix{X\ar[r]^{\phi}\ar[d]^{\id}&[r_{n-1}]\ar[r]^{\id}\ar[d]_{\tau}&
  [r_{n-1}]\ar[r]^{f_{n-1}} \ar[d]^{\sigma_{n-1}}&
\ldots\ar[r]^{f_2}&[r_1]\ar[d]^{\sigma_1}\\
X\ar[r]^{\tau\phi}&[r_{n-1}]\ar[r]^{g_n}& [s_{n-1}]\ar[r]^{g_{n-1}}&
\ldots\ar[r]^{g_2}&[s_1].}$$ Let  $(x_0,\ldots,x_{s_{n-1}})$ be the
sequence of cardinalities of the  preimages of $\sigma_{n-1}$, which
also determines $g_n$. There exists a  bijection $\xi$ of
$[r_{n-1}]$ such that $\sigma_{n-1}=g_n\xi$. If $\xi,\xi'$ are
bijections of $[r_{n-1}]$ both satisfying  the previous  equality
then
$\xi(\xi')^{-1}\in\Sigma_{x_0}\times\ldots\times\Sigma_{x_{s_{n-1}}}$.
Any element $\tau$ that makes the diagram commute is of the form
$\alpha\xi$  for
$\alpha\in\Sigma_{x_0}\times\ldots\times\Sigma_{x_{s_{n-1}}}$. As in
proposition \ref{prop:shufflecycles}, the element
$\mathrm{sgn}(\xi;X)(c_{(x_0,\ldots,x_{s_{n-1}});X})\xi$ does not
depend on the choice  of $\xi$ and it is a generator of
$H_{(r_n,s_{n-1},\ldots,s_1)}({\epi_n^X}^{t,\phi},\partial_n)$. This
gives the desired isomorphism of $k$-modules between this homology
group and
$C_{(s_{n-1},\ldots,s_1)}^{E_{n-1},X[1]}({\epi_{n-1}^{X[1]}}^{\bar
t,\phi})$.

A direct inspection of the signs in
\ref{def:enchainsX} shows that the induced differential $\partial_i$
coincides with the one on
$C_{(s_{n-1},\ldots,s_1)}^{E_{n-1},X[1]}({\epi_{n-1}^{X[1]}}^{\bar
t,\phi})$ for $1\leq i\leq n-2$. The case $i=n-1$ is similar to
the proof of corollary \ref{cor:fork2}: let us choose a
bijection $\xi\colon [r_{n-1}]\rightarrow [r_{n-1}]$ with the property
that $i<j$ and   
$g_n\xi(i)=g_n\xi(j)$ imply $\xi(i)<\xi(j)$. We take
$\mathrm{sgn}(\xi;X)(c_{(x_0,\ldots,x_{s_{n-1}});X})\xi$ as the
corresponding generator.

On the one hand, one has
$\partial_{n-1}(\xi,\sigma_{n-1},\ldots,\sigma_1)=\sum_{j|g_{n-1}(j)=g_{n-1}(j+1)}
(-1)^{s_{n-1,j}}
\epsilon_j\alpha_j$ where $\alpha_j$ is a sum of trees which is a
generator of the form
$\mathrm{sgn}(\xi;X)c_{(x_0,\ldots,x_j+x_{j+1},\ldots,x_{s_{n-1}});X}\xi$ 
and $\epsilon_j$ is a sign to be determined. If $Y_j=\{y_1<\ldots<
y_{x_j}\}$ and $Y_{j+1}=\{z_1<\ldots <z_{x_{j+1}}\}$ are the sets
$(g_{n-1}g_n\xi\phi)^{-1}(j)$ and $(g_{n-1}g_n\xi\phi)^{-1}(j+1)$
respectively then the sign $\epsilon_j$ is precisely
$\epsilon(Y_j[1];Y_{j+1}[1])$. On the other hand, computing
$$\partial_{n-1}\left(\vcenter{\xymatrix{X[1]\ar[r]^{\phi}\ar[d]^{\id}&
  [r_{n-1}]\ar[r]^{f_{n-1}} \ar[d]^{\sigma_{n-1}}&
\ldots\ar[r]^{f_2}&[r_1]\ar[d]^{\sigma_1}\\
X[1]\ar[r]^{\sigma_{n-1}\phi}& [s_{n-1}]\ar[r]^{g_{n-1}}&
\ldots\ar[r]^{g_2}&[s_1]}}\right)$$ gives the same signs for
$Y_j[1]=(g_{n-1}\sigma_{n-1}\phi)^{-1}(j)=(g_{n-1}g_n\xi\phi)^{-1}(j)$.

\end{proof}

\end{document}